\documentclass{elsart}	
\usepackage{amsmath,  amssymb, amscd}
\usepackage[all]{xy}
\journal{{\bf Advances in Mathematics} (accepted on 29/10/04)}
\newcommand{\A}{\mathcal{A}}
\newcommand{\B}{\mathcal{B}}

\newcommand{\K}{\mathcal{K}}
\newcommand{\cL}{\mathcal{L}}
\renewcommand{\a}{\alpha}
\renewcommand{\b}{\beta}

\renewcommand{\d}{\delta}

\newcommand{\tr}{\text{tr}\,}

\newcommand{\extd}{\text{d}}

\newcommand{\ad}{\text{ad}}
\newcommand{\Ch}{\mathrm{Ch}}

\theoremstyle{plain}
\newtheorem{theorem}{Theorem}[section]
\newtheorem{lemma}[theorem]{Lemma}
\newtheorem{proposition}[theorem]{Proposition}
\newtheorem{corollary}[theorem]{Corollary}
\theoremstyle{remark}

\newtheorem{remark}[theorem]{Remark}
\theoremstyle{definition}
\newtheorem{definition}{Definition}
\begin{document}
\date{}
\begin{frontmatter}

\title{On a generalized Connes-Hochschild-Kostant-Rosenberg theorem}

\author[Adelaide]{Varghese Mathai\thanksref{ARC}}
\ead{vmathai@maths.adelaide.edu.au }
and
\author[Adelaide]{Danny Stevenson\thanksref{UCR}\thanksref{ARC}}
\ead{dstevens@math.ucr.edu }

\address[Adelaide]
{Department of Pure Mathematics,
University of Adelaide,
Adelaide, SA 5005,
Australia}

\thanks[UCR]
{Currently: Department of Mathematics,
University of California, Riverside, CA}

\thanks[ARC]{The authors acknowledge the support of the Australian
Research Council.}

\begin{abstract}
The central result of this paper is an explicit computation of the 
Hochschild and
cyclic homologies of a natural smooth
subalgebra of stable continuous trace algebras
having smooth manifolds $X$ as their spectrum.
More precisely, the Hochschild homology
is identified with the space of differential forms on $X$, and the
periodic cyclic homology with the twisted de Rham cohomology of $X$,
thereby generalizing some fundamental
results of Connes and Hochschild-Kostant-Rosenberg. The
Connes-Chern character is also identified here with the twisted Chern
character.
\end{abstract}

\end{frontmatter}
%\maketitle

%-----------------------------------------------------------------------

\section{Introduction}

In \cite{Connes}, Connes showed how to extend both classical
Hochschild homology and  cyclic homology to the
category of locally convex algebras,
computing in particular the Hochschild homology
of the Fr\'{e}chet algebra of smooth functions on a
compact manifold $X$.  This he showed was canonically isomorphic to the
space of differential forms
on $X$, thus generalizing a fundamental result of
Hochschild-Kostant-Rosenberg \cite{HKR}.
In the same paper, Connes also
identified the periodic cyclic homology
of the Fr\'{e}chet algebra of smooth functions on $X$
with the de Rham cohomology of $X$.  Connes's
theorem shows that cyclic homology is indeed a far-reaching
generalisation of de Rham cohomology for noncommutative
topological algebras, and serves as a
central motivation for the results here.

Since \cite{Connes}, there have been extensions
in different directions including
the case of non-compact manifolds due to Pflaum \cite{Pfl}, the 
equivariant
case studied by Block and Getzler \cite{BG} and the case of \'etale 
groupoids and foliations
\cite{BN,C}.
In this paper we take up the suggestion of \cite{BM} and
prove the following generalisation of the
Connes and Hochschild-Kostant-Rosenberg theorem.
Let $A$ be a stable continuous trace $C^*$-algebra with spectrum
a smooth, compact manifold $X$.
Then by \cite{DD},
$A =C(X, \K(P) )$ is the algebra of continuous sections of a
locally
trivial bundle $\K(P) = P \times_{PU} \K$ on $X$ with fibre the
algebra $\K$ of compact operators on a separable
Hilbert space associated to a principal $PU$
bundle $P$ on $X$ via the adjoint action of
$PU$ on $\K$.  Here $PU$ is the group
of projective unitary operators on the Hilbert space.
Such algebras  $A$
are classified up to isomorphism by their
Dixmier-Douady invariant   $\delta(P)\in H^3(X;\Zset)$.
Standard approximation theorems show that the transition
functions of $P$ can be taken to be smooth and so inside $A$ we can 
consider
a dense, canonical smooth $*$-subalgebra $\A = C^\infty(X, \cL^1(P))$
consisting of all
smooth sections of the sub-bundle $\cL^1(P) = P \times_{PU} \cL^1$ of
$\K(P)$ with
fibre the algebra  $\cL^1$ of trace class operators on the Hilbert space
and structure group $PU$.  Then a central result in this
paper is the following,

\begin{theorem}\label{main}
The Hochschild homology $HH_\bullet(\A)$ of
the Fr\'{e}chet algebra $\A$ is isomorphic to the
space of differential forms $\Omega^\bullet(X)$ and
the periodic cyclic homology $HP_\bullet(\A)$ is isomorphic to the
twisted de Rham cohomology $H^\bullet(X;c(P))$
for some closed $3$-form $c(P)$ on $X$ such that
$\frac{1}{2\pi i}c(P)$ represents the image of $\delta(P)$
in real cohomology.
\end{theorem}

Recall from \cite{BCMMS} that
the twisted de Rham cohomology $H^\bullet(X;c(P))$ of the
manifold $X$ is defined as the cohomology of the complex
$\left(\Omega^\bullet(X), \extd - c\right)$, where $c$ denotes exterior
multiplication
by the closed 3-form $c = c(P)$.
By choosing a connection $\nabla$ on $\cL^1(P)$ satisfying a
derivation property with respect to the algebra structure
we construct a chain map $\Ch$ from the
complex computing periodic cyclic homology to the
twisted de Rham complex $(\Omega^\bullet(X),
\extd -c)$.  By examining the behaviour of this
map when we perturb the connection $\nabla$ we
are able to deduce that this chain map is locally
a quasi-isomorphism on $X$.  Standard double complex
arguments then give that $\Ch$ is a quasi-isomorphism
globally.   Extensions to
$*$-algebras of smooth sections of a smooth algebra bundle
with fibre the algebra of operators belonging to the
Schatten class ideal $\cL^p$
are also discussed.

Recall from \cite{Ros} that  the twisted $K$-theory  is by
definition $K^i(X, P) = K_i(A),$ which is in turn
isomorphic to $K_i(\A)$ as argued in section 4.2. In
\cite{BCMMS,MS},  it was established that the twisted $K$-theory
$K_\bullet (\A)$ could be identified with
so-called `bundle gerbe
$K$-theory', $K_{bg}^\bullet(X,P)$ (also denoted by $K^\bullet(X,P)$),
thus giving a geometric description of
elements of twisted $K$-theory. Using this geometric description, it was
established that there is a twisted Chern character $\text{ch}_P :
K^\bullet(X, P) \to
H^\bullet(X,c(P))$ possessing certain functorial properties.

It turns out
that the Connes-Chern character homomorphism
      $\text{ch}\colon K_i(\A)\to HP_i(\A)$ satisfies the same
functorial properties.  This is not surprising, once we
establish that
under the natural identifications $K_\bullet (\A)$ with
$K^\bullet(X,P)$, and
$HP_\bullet(\A)$ with $H^\bullet(X,c(P))$, the Connes-Chern character
homomorphism
agrees with the twisted Chern character homomorphism
constructed in \cite{BCMMS}, which
can be encapsulated
as the commutativity of the following diagram,
\begin{equation} \label{main diag}
\begin{CD}
K^\bullet(X,  P)  @>\cong >> K_\bullet (\A) \\
           @V{\text{ch}_P}VV          @VV{\text{ch}} V     \\
H^\bullet(X, c(P))   @>\cong >>  HP_{\bullet} (\A) .
\end{CD}\end{equation}

We now outline the paper. In
section $2$ we first give some
standard background material on
topological algebras and establish a result
(Proposition~2.3) that we will need later
in the paper.  In section 3 we review
preliminary material on defining the Hochschild, cyclic
and periodic cyclic homology groups $HH_\bullet(\A)$,
$HC(\A)$ and $HP_\bullet(\A)$ respectively for a Fr\'{e}chet algebra
$\A$, and the computation of some examples, such as
$HH_\bullet(\cL^1),
HC_\bullet(\cL^1)$ and $HP(\cL^1)$.
In Section 4, we review the $K$-theory of such algebras, and specialize
to the case of interest
for the rest of the paper, namely $\A = C^\infty(X, \cL^1(P))$, where
$P$ is a principal $PU$-bundle
with Dixmier-Douady class equal to $[c] \in H^3(X, \Zset)$, and $
\cL^1(P)
= P \times_{PU} \cL^1$
is the bundle with fibre $\cL^1$ associated to
$P$ via the adjoint action of $PU$ on $\cL^1$.
We show that the $K$-theory $K_\bullet(\A)$ can be
identified with the twisted $K$-theory $K^\bullet(X;P)$.
We also review the construction of the
Connes-Chern character $\text{ch}: K_\bullet(\A) \to
HP_\bullet(\A)$.
In Section 5 we construct, following ideas of
Gorokhovsky \cite{Gor1,Gor2}, a chain map
$\mathrm{Ch}_\bullet\colon CC_\bullet(\A)\otimes \Cset((u))
\to \Omega^\bullet(X)\otimes \Cset((u))$ from the
complex $CC_\bullet(\A)\otimes \Cset((u))$
computing periodic cyclic homology to
the twisted de Rham complex $(\Omega^\bullet(X)\otimes
\Cset((u)),\extd + uc)$.
We remark that the map $\mathrm{Ch}$ can be
constructed in the formalism of
Quillen \cite{Qui}.
In Theorem~\ref{thm: main result}
we establish that the map $\mathrm{Ch}_\bullet$
is a quasi-isomorphism.  We also describe extensions
of this result to Hochschild and cyclic homology
respectively.
The arguments we use here are
reminiscent of the \v{C}ech-de Rham tic-tac-toe type argument.
Finally in Section 6, we prove the commutativity of the diagram
\eqref{main diag}.  We also prove that the
Connes-Chern map $\text{ch}\colon K_i(\A)\to HP_i(\A)$
becomes an isomorphism on tensoring with the
complex numbers.  \\

\noindent\textbf{Acknowledgments.} It is a pleasure for
the authors to thank Joachim Cuntz for his comments and
very useful references and Ezra Getzler for his comments
on reading a first draft of the paper and
suggesting the elegant rephrasing of several results in the paper.

%-----------------------------------------------------------------------

\section{Preliminaries on Topological Algebras}

A \emph{Fr\'{e}chet algebra} $A$ is an algebra with a complete
locally convex topology
determined by a
countable family of semi-norms $\{p_{\a}\}$
which are required to be sub-multiplicative
in the sense that $p_{\a}(ab)\leq p_{\a}(a)
p_{\a}(b)$.  This implies that the multiplication
map $A\times A\to A$ is jointly continuous.
We could in fact consider a larger class of
algebras, so-called $m$-algebras for which we refer
to Cuntz's article \cite{Cu1}, but for the purposes of
this note Fr\'{e}chet algebras will suffice.

Let us also remark that the unitization $\tilde{A}$ of
a Fr\'{e}chet algebra $A$ is again a Fr\'{e}chet algebra.
Recall that $\tilde{A} = A\oplus \Cset$ with multiplication
given by $(a,\lambda)\cdot (b,\mu) = (ab+ \lambda b+
\mu a,\lambda\mu)$.  The topology of $\tilde{A}$ is determined
by the family of semi-norms $\{\tilde{p}_{\a}\}$ given by
$\tilde{p}_{\a}(a,\lambda) = p_{\a}(a) + |\lambda|$.  Similarly
the direct sum $A\oplus B$ of two Fr\'{e}chet algebras
is again Fr\'{e}chet.  One needs to take
more care in defining tensor products, the tensor product
in which we will be exclusively interested in is
Grothendieck's
projective tensor product $\otimes_{\pi}$ \cite{Tre}.  Recall that this
is defined as follows: if $A$ and $B$ are Fr\'{e}chet
algebras with countable families of semi-norms $\{p_{\a}\}$
and $\{q_{\b}\}$ then the algebraic tensor product $A\odot B$ is
given a locally convex topology, called the \emph{projective tensor
product}
topology, defined by the countable family of semi-norms
given by
$$
p_{\a}\otimes q_{\b}(c) = \inf\{\sum^n_{i=1}p_{\a}(a_i)q_{\b}(b_i)|\
c=\sum^n_{i=1}a_i\odot b_i\}
$$
We will denote the algebraic tensor product $A\odot B$
equipped with this topology by $A\otimes_{\pi}B$ and the
completion of $A\otimes_{\pi}B$ by $A\hat{\otimes}_{\pi}B$.
$A\hat{\otimes}_{\pi}B$ is again a Fr\'{e}chet algebra.

Recall the well known example of the Fr\'{e}chet algebra $C^\infty(X)$
of smooth functions on a Riemannian manifold $X$.
The seminorms on this algebra are defined as follows.
Let $V$ be a chart in $X$ and $\{K_n\}_{n\in \mathbb N}$
be an exhausting sequence of compact subsets
of $V$, i.e., $K_n$ is contained in the interior of $K_{n+1}$
and every compact subset of $V$ is contained in some
$K_n$. For $f \in C^\infty(X)$, we denote its restriction to $V$
(and to $K_n$) by the same symbol. Then set
%This can be equipped
%with the (countable) family of  semi-norms
%$\{p_{n, V}\}$ where $n\in \mathbb N$, which are defined using local
%charts $V$ by
$$
p_{n, V}(f) = \sum_{j=0}^n \frac{1}{j!} \sup_{x \in K_n}
|\partial^\alpha f(x)|.
$$
where the partial derivatives are taken in the chart $V$.
Then varying over all charts $V$ and all $n \in \mathbb N$, $p_{n, V}$
forms a basis of continuous seminorms for $C^\infty(X)$,
     and the derivation property of
the partial derivatives ensures that the family of semi-norms
is sub-multiplicative, giving it the structure of a Fr\'{e}chet algebra.

Any Banach algebra $\A$ is a Fr\'{e}chet algebra with
the family of semi-norms consisting of the single norm
$|\cdot|$.  In particular we will be interested in the
case where $\A = \cL^1(H)$, the Banach algebra of trace
class operators on a separable Hilbert space $H$.  We
take this opportunity to review some facts about $\cL^1(H)$.
Let $H$ be an infinite dimensional Hilbert space and $\{v_n\}, \;
n\in \Nset$
an orthonormal basis of vectors in $H$. Let $\B = \B(H)$ denote the
bounded
operators
on $H$. We can define a semi-finite trace on positive elements in
$\B$
by
\begin{equation}\label{eq:tr}
\tr(x) = \sum_{n=1}^\infty (x, v_n).
\end{equation}
It turns out that $\tr(x)$ is well defined, independent of the choice
of orthonormal
basis and satisfies $\tr(x^*x) = \tr(xx^*)$ for all $x \in \B$.
For $p\ge 1$, define the Schatten class
$$
\cL^p = \left\{ x \in \B: \tr(|x|^p) < \infty\right\}.
$$
where $|x| = \sqrt{x^*x}$. Then it is well known that $\cL^p$ is a ideal
in $\B$ that
is closed under taking adjoints, and which is contained in the closed
ideal of compact
operators $\K$ in $\B$. $\cL^p$ is a Banach $*$-algebra with respect to
the norm
$||x||_p = (\tr(|x|^p))^{1/p}$, but it has no unit.
$\cL^1$ is the Banach $*$-algebra of trace class operators and
$\cL^2$ the Banach $*$-algebra of Hilbert-Schmidt operators.
Moreover, it is the case that the trace map is well defined for
arbitrary
elements of $\cL^1$,
\begin{equation}\label{eq:tr2}
\tr : \cL^1 \to \Cset
\end{equation}
and is given by the same formula \eqref{eq:tr}.
Moreover $\cL^p \subset \cL^q$ whenever $p \le q$.

\begin{lemma}\label{lem:tensor}
Suppose that $X$ and $Y$ are paracompact manifolds.   Then one has
a canonical isomorphism,
$$
C^\infty(X,\cL^1)\hat{\otimes}_{\pi}C^\infty(Y,\cL^1) \cong
C^\infty(X\times Y,\cL^1\,\hat{\otimes}_{\pi}\,\cL^1).
$$
\end{lemma}

\begin{pf*}{Proof.}
Recall that the projective tensor product satisfies
the following two properties (see for example \cite{Tre}).
If $X$ and $Y$ are smooth manifolds then $C^\infty(X)
\hat{\otimes}_{\pi}C^\infty(Y) = C^\infty(X\times Y)$
and if $E$ is a Fr\'{e}chet space then $C^\infty(X)
\hat{\otimes}_{\pi}E = C^\infty(X,E)$, which is the
Fr\'{e}chet space of $E$-valued smooth functions on $X$.

Therefore $C^\infty(X,\cL^1) \cong C^\infty(X) \hat
\otimes_\pi \cL^1$
and $C^\infty(Y,\cL^1) \cong C^\infty(Y) \hat \otimes_\pi \cL^1$, and
$$
\begin{array}{lcl}
C^\infty(X,\cL^1)\hat{\otimes}_{\pi}C^\infty(Y,\cL^1) &\cong&
C^\infty(X) \hat \otimes_\pi    C^\infty(Y) \hat \otimes_\pi \cL^1
\hat{\otimes}_{\pi}        \cL^1\\
& \cong & C^\infty(X\times Y)  \hat \otimes_\pi \cL^1
\hat{\otimes}_{\pi}       \cL^1\\
& \cong & C^\infty(X\times Y,\cL^1\,\hat{\otimes}_{\pi}\,\cL^1).
\end{array}
$$
\end{pf*}

\begin{remark}
We caution the reader that the projective tensor product,
$\cL^1\,\hat{\otimes}_{\pi}\,\cL^1$
is not  isomorphic to $\cL^1$. A similar remark applies when $\cL^1$ is
replaced by $\cL^p$ or
by $\K$. This is in contrast to the well known fact that
$\K\,{\otimes}\,\K
\cong \K$, when the $C^*$ tensor product is applied instead.
\end{remark}

Another Fr\'{e}chet algebra that we will
be interested in is $C^\infty(X,\cL^1(P))$, the
algebra of smooth
sections of a  Banach algebra bundle $\cL^1(P) = P\times_{PU} \cL^1$
with fibre the algebra of trace class operators $\cL^1$
associated to a principal $PU$ bundle $P$ via the
adjoint action of $PU$ on $\cL^1$.
Here $PU = PU(H)$ is the projective unitary group
of $H$ equipped with the norm topology (see \S 4.2 for more details).
      The trace $\tr(\cdot)\colon \cL^1(H)\to
\mathbb{C}$ can be extended fibre-wise
$\tr(\cdot)\colon C^\infty(\cL^1(P))\to C^\infty(X)$,
and the algebra
of sections $C^\infty(\cL^1(P)) = C^\infty(X,\cL^1(P))$
is made into a
Fr\'{e}chet algebra as follows.
Let $V$ be a chart in $X$ and $\{K_n\}_{n\in \mathbb N}$
be an exhausting sequence of compact subsets
of $V$. For $f \in C^\infty(X,\cL^1(P))$, we denote its restriction to
$V$
(and to $K_n$) by the same symbol. Then set
$$
p_{n, V}(f) = \sum_{j=0}^n \frac{1}{j!} \sup_{x \in K_n} \tr(|\nabla^j
f(x)|)
$$
where $\nabla^j f$ is the $j$-th covariant derivative
of the section $f$  taken the the chart $V$, $\nabla$ is a connection
on the bundle $\cL^1(P)$ which acts as a derivation
on $C^\infty(X,\cL^1(P))$ in the sense
that
$\nabla(fh) = f \nabla h + \nabla f h$ (the existence
of such connections will be established in \S 5.1.).
     Then varying over all charts $V$ and all $n \in \mathbb N$, $p_{n,
V}$
forms a basis of continuous seminorms for $C^\infty(X)$,
giving it the structure of a Fr\'{e}chet algebra.
%Here
%$||\nabla^j f(x)||_1 = \tr(|\nabla^j f(x)|)$ is the pointwise
%trace norm.
If $\nabla'$ is another connection on $\cL^1(P)$ that is also
compatible with the product as above, then $\nabla'$ defines a
quasi-isometric family of semi-norms $\{p'_{n, V}\}$ provided that
$\nabla$ and $\nabla'$ differ by a $1$-form on $X$ which takes
values in the sub-bundle $\B(\cL^1(P))\subseteq
\mathrm{End}(\cL^1(P))$ of \emph{bounded}
endomorphisms of $\cL^1(P)$.  The derivation property of
$\nabla$ ensures that the family of semi-norms
is sub-multiplicative.

The following is a generalization of Lemma \ref{lem:tensor} to the case
of non-trivial bundles of algebras with fibre $\cL^1$ and structure
group $PU$.

\begin{proposition}\label{prop:tensor}
Suppose that $P$ and $Q$ are principal $PU$ bundles
on paracompact manifolds $X$ and $Y$ respectively.  Denote
by $\cL^1(P)$ and $\cL^1(Q)$ the Banach vector bundles on
$X$ and $Y$ respectively with fibre $\cL^1$ associated to
$P$ and $Q$ via the adjoint action of $PU$ on $\cL^1$.
Suppose also $\cL^1(P)$
and $\cL^1(Q)$ come
equipped with connections $\nabla_{\cL^1(P)}$ and $\nabla_{\cL^1(Q)}$
satisfying the derivation property above.
If $C^\infty(X,\cL^1(P))$ and $C^\infty(Y,\cL^1(Q))$ are equipped with
the topologies coming from the families of semi-norms
induced by $\nabla_{\cL^1(P)}$ and $\nabla_{\cL^1(Q)}$, then we have
$$
C^\infty(X,\cL^1(P))\hat{\otimes}_{\pi}C^\infty(Y,\cL^1(Q)) =
C^\infty(X\times Y,\cL^1(P)\,\boxtimes\,\cL^1(Q)),
$$
where $\cL^1(P)\,\boxtimes\,\cL^1(Q)$ is the algebra bundle over the
cartesian product $X \times Y$,
with fibre $\left(\cL^1(P)\,\boxtimes\,\cL^1(Q)\right)_{(x, y)}
= \cL^1(P)_x \hat{\otimes}_{\pi}\,\cL^1(Q)_y$,
i.e. modelled on $\cL^1\,\hat{\otimes}_{\pi}\,\cL^1$.
\end{proposition}

\newcommand{\cE}{\mathcal{E}}
\newcommand{\cF}{\mathcal{F}}
\begin{pf*}{Proof.}
We note first that $C^\infty(\cL^1(P))\odot C^\infty(\cL^1(Q))$
can be identified with a subspace of $C^\infty(\cL^1(P)\boxtimes
\cL^1(Q))$.
We will prove two things:
\begin{enumerate}
\item $C^\infty(\cL^1(P))\otimes_{\pi}
C^\infty(\cL^1(Q))\subset C^\infty(\cL^1(P)\boxtimes \cL^1(Q))$ is a
homeomorphism onto its image,

\item $C^\infty(\cL^1(P))\odot
C^\infty(\cL^1(Q))$ is dense in $C^\infty(\cL^1(P)\boxtimes \cL^1(Q))$.
\end{enumerate}
This suffices to show the equality
$C^\infty(\cL^1(P))\hat{\otimes}_{\pi}
C^\infty(\cL^1(Q)) = C^\infty(\cL^1(P)\boxtimes \cL^1(Q))$.
To show that
$C^\infty(\cL^1(P))\otimes_{\pi}C^\infty(\cL^1(Q))\subset
C^\infty(\cL^1(P)
\boxtimes \cL^1(Q))$ is a homeomorphism onto its image, first of all
observe that the inclusion
$C^\infty(\cL^1(P))\otimes_{\pi}C^\infty(\cL^1(Q))\subset
C^\infty(\cL^1(P)\boxtimes \cL^1(Q))$ is continuous.  Now suppose
we have a sequence $f_i$ in $C^\infty(\cL^1(P))\odot
C^\infty(\cL^1(Q))$ which converges to zero in the
$C^\infty(\cL^1(P)\boxtimes \cL^1(Q))$ topology.  We have
to show that $f_i \to 0$ in the projective tensor
product topology on $C^\infty(\cL^1(P))\odot
C^\infty(\cL^1(Q))$.  Choose good open
covers $\{U_{\a}\}_{\a\in A}$ and $\{V_{\b}\}_{\b\in B}$
of $X$ and $Y$ respectively and denote by $\cL^1(P)_{\a}$ the
restriction
$\cL^1(P)|_{U_{\a}}$ and by $\cL^1(Q)_{\b}$ the restriction of
$\cL^1(Q)|_{V_{\b}}$.  Since $U_{\a}, V_{\b}$ are contractible, the
restricted bundles $\cL^1(P)_{\a}, \cL^1(Q)_{\b}$ are trivializable,
therefore
$C^\infty(\cL^1(P)_{\a}) \cong C^\infty(U_\a, \cL^1)$ and
$C^\infty(\cL^1(Q)_{\b})  \cong C^\infty(V_\b, \cL^1)$.
By remarks above and by Lemma \ref{lem:tensor}, we have the canonical
isomorphism
$$C^\infty(\cL^1(P)_{\a})\hat{\otimes}_{\pi}
C^\infty(\cL^1(Q)_{\b})  = C^\infty(\cL^1(P)_{\a}
\boxtimes \cL^1(Q)_{\b}).$$  Denote by
$r_{\a}\colon C^\infty(\cL^1(P))\to
C^\infty(\cL^1(P)_{\a})$, $r_{\b}\colon C^\infty(\cL^1(Q))
\to C^\infty(\cL^1(Q)_{\b})$ and $r_{\a,\b}\colon
C^\infty(\cL^1(P)\boxtimes \cL^1(Q))\to
C^\infty(\cL^1(P)_{\a}\boxtimes \cL^1(Q)_{\b})$
the restriction maps.  Then we have a commutative diagram
\begin{equation}
\label{eq: tensor comm diagram}
\xymatrix{
C^\infty(\cL^1(P)){\otimes}_{\pi}
C^\infty(\cL^1(Q)) \ar[d]_-{\rm inc } \ar[r]^-{\mathrm{inc}} &
C^\infty(\cL^1(P)\boxtimes \cL^1(Q)) \ar[d]^-{||} \\
C^\infty(\cL^1(P))\hat{\otimes}_{\pi}
C^\infty(\cL^1(Q)) \ar[d]_-{\prod (r_{\a}
\hat{\otimes}_{\pi}r_{\b})} &
C^\infty(\cL^1(P)\boxtimes \cL^1(Q)) \ar[d]^-{r_{\a,\b}} \\
\displaystyle{\prod_{\a,\b}}
C^\infty(\cL^1(P)_{\a})\hat{\otimes}_{\pi}
C^\infty(\cL^1(Q)_{\b}) \ar[r]^-{=} &
\displaystyle{\prod_{\a,\b}}
C^\infty(\cL^1(P)_{\a}\boxtimes
\cL^1(Q)_{\b})                          }
\end{equation}
where $\prod$ denotes the direct product.
Both the lower vertical maps are split injective (a
splitting can be constructed via partitions
of unity subordinate to the open covers $\{U_{\a}\}_{\a\in A}$
and $\{V_{\b}\}_{\b\in B}$).
%Here we regard the completion
%$C^\infty(\cL^1(P))\hat{\otimes}_{\pi}C^\infty(\cL^1(Q))$
%as sitting inside $C^\infty(\cL^1(P)\boxtimes \cL^1(Q))$.
%\marginpar{\tiny{Do you see any problem with this Mathai,
%can we regard $C^\infty(\cL^1(P))\hat{\otimes}_{\pi}
%C^\infty(\cL^1(Q))$ as sitting inside $C^\infty(\cL^1(P)
%\boxtimes \cL^1(Q))$?}}
It follows from the commutativity of
the diagram~\eqref{eq: tensor comm diagram}
that $r_{\a}\hat{\otimes}_{\pi}r_{\b}(f_i)\to 0$
and hence $f_i\to 0$ since $\prod r_{\a}\hat{\otimes}_{\pi} r_{\b}$
is split injective.  Therefore the inclusion
$C^\infty(\cL^1(P))\otimes_{\pi}C^\infty(\cL^1(Q))
\subset C^\infty(\cL^1(P)\boxtimes \cL^1(Q))$ is a
homeomorphism onto its image.
We now have to show that $C^\infty(\cL^1(P))\odot
C^\infty(\cL^1(Q))$ is dense in $C^\infty(\cL^1(P)
\boxtimes \cL^1(Q))$.
Let $f\in C^\infty(\cL^1(P)
\boxtimes \cL^1(Q))$ and let $f_{\a,\b} = r_{\a,\b}(f)$
denote its image in $C^\infty(\cL^1(P)_{\a}
\boxtimes \cL^1(Q)_{\b}) = C^\infty(\cL^1(P)_{\a})
\hat{\otimes}_{\pi}C^\infty(\cL^1(Q)_{\b})$
under the restriction map
$r_{\a,\b}$.  Choose a sequence $f^i_{\a,\b}\in
C^\infty(\L^1(P)_{\a})\odot C^\infty(\cL^1(Q)_{\b})$
converging to $f_{\a,\b}$.  Do this for all $\a$
and $\b$ and denote by $f_i$ the image of the
sequence $\prod f^i_{\a,\b}$ under the
splitting for the map $\prod r_{\a}\odot r_{\b}$.  It
follows from the commutativity of the
diagram~\eqref{eq: tensor comm diagram} above
that $f_i$ converges to $f$ in $C^\infty(\cL^1(P)\boxtimes
\cL^1(Q))$.
\end{pf*}

%-----------------------------------------------------------------------

\section{Cyclic Homology}
\label{sec:one}
\subsection{Definitions}
\label{sec:one a}

In this section we recall the definition and
main properties of Hochschild, cyclic and
periodic cyclic homology.  We shall closely follow
the articles of Block and Getzler \cite{BG} and
Cuntz \cite{Cu1}.
Let $A$ be a unital algebra over the complex numbers $\Cset$
and let $CC_k(A) = A^{\otimes k+1}$.  Define
a differential $b$ of degree $-1$ by
\begin{multline}
b(a_0\otimes a_1\otimes \cdots \otimes a_k) = \\
a_0a_1\otimes a_2\otimes \cdots \otimes a_k +
\sum^n_{j=2}(-1)^{j-1}a_0\otimes \cdots \otimes
a_{j-1}a_j\otimes \cdots \otimes a_k \\
+ (-1)^n a_na_0\otimes a_1\otimes
\cdots \otimes a_{k-1}.  \label{eq: formula for Hoch b}
\end{multline}
$b$ is defined to be zero on $CC_0$.  The \emph{Hochschild
homology} $HH_\bullet(A)$ is the homology of the
complex $(CC_\bullet(A),b)$.
This definition extends to the non-unital case
as described in Loday \cite{Lod} pages 29--30 by
setting
$$HH_n(A) =  \text{coker} \left(HH_n(\Cset) \to
HH_n(\tilde{A})\right)$$
for a non-unital algebra $A$, where $\tilde{A}$ denotes the
unitization of $A$.  If however, $A$ is
\emph{$H$-unital}, then $HH_n(A)$ can be computed
from the same complex $(CC_\bullet(A),b)$ which
computes Hochschild homology in the unital case.
Recall \cite{Lod} that $A$ is said to be $H$-unital
if the bar complex $\mathrm{Bar}_\bullet(A)$ is
contractible.  The bar complex
has $\mathrm{Bar}_k(A) = A^{\otimes k}$ in degree
$k$ with differential $b'$ given by the standard
formula
$$
b'(a_0\otimes a_1\otimes \cdots \otimes a_{k-1}) =
\sum^{k-2}_{j=0}(-1)^j a_0\otimes \cdots \otimes a_ja_{j+1}
\otimes \cdots \otimes a_{k-1}.
$$
Another way to compute
$HH_n(A)$, which works regardless of whether $A$
is unital or not, is to use the \emph{reduced}
Hochschild complex which is defined to be
$\overline{CC}_k(A) = \tilde{A}\otimes A^{\otimes k} =
A^{\otimes k}\oplus A^{\otimes k+1}$ in
degrees $k>0$ and $\overline{CC}_0(A) = A$ with differential
$b$ given by the obvious extension
of the formula~\eqref{eq: formula for Hoch b} above.
More specifically, the operator $b\colon \overline{CC}_{k+1}(A)
\to \overline{CC}_k(A)$ can be written \cite{Cu1} with
respect to the splitting $\overline{CC}_{k+1}(A) =
A^{\otimes k+1} \oplus A^{\otimes k+2}$ as
$$
b =\left( \begin{array}{cc}
b & 1- \lambda \\
0 & - b'
\end{array} \right),
$$
where $\lambda\colon A^{\otimes k}\to A^{\otimes k}$
is Connes' signed permutation operator
$\lambda(a_0\otimes a_1\otimes \cdots\otimes a_{k-1})
= (-1)^k a_{k-1}\otimes a_0\otimes \cdots \otimes a_{k-2}$.
       From now on we will write $CC_k(A)$ for
$\overline{CC}_k(A)$.
We also introduce Connes' operator $B$ of degree $1$ defined by
$$
B(\tilde{a}_0\otimes a_1\otimes \cdots \otimes a_k) =
\sum^k_{i=0}(-1)^{ik} 1\otimes a_i \otimes a_k\otimes
a_0\otimes \cdots \otimes a_{i-1}.
$$
$b$ and $B$ satisfy the identities $b^2 = B^2 = bB + Bb = 0$.
We can also write $B$, with respect to the splitting
$CC_k(A) = A^{\otimes k+1}\oplus A^{\otimes k}$, as
$$
B = \left( \begin{array}{cc}
0 & 0 \\
Q & 0
\end{array} \right)
$$
where $Q = \sum^k_{j=0} \lambda^k$.
Following Block and Getzler, we introduce a variable
$u$ of degree $-2$.  We then define the \emph{cyclic homology}
$HC_\bullet(A)$ of $A$ to be the homology
of the complex
$$
(CC_\bullet(A)\otimes \Cset((u))/u\Cset[[u]], b+uB)
$$
while the \emph{periodic cyclic homology} $HP_\bullet(A)$
of $A$ is the homology of the complex
$$
(CC_\bullet(A)\otimes \Cset((u)),b+uB).
$$
Note that the periodic theory
$HP_\bullet(A)$ is a $\Zset$-graded version of the
usual $\Zset/2$-graded version.
In this note we will be concerned with Hochschild, cyclic and
periodic cyclic homology defined for \emph{topological} algebras.
Following \cite{Connes} we shall define \emph{topological}
Hochschild, cyclic and periodic cyclic homology for
Fr\'{e}chet algebras $\A$.
Suppose that $\A$ is a Fr\'{e}chet
algebra.
Then we form the vector spaces $CC_k(\A) =
\tilde{\A}\hat{\otimes}_{\pi}\A^{\hat{\otimes}_{\pi} k}$
by replacing the ordinary tensor products
$\otimes$ with the completed projective
tensor product $\hat{\otimes}_{\pi}$
of Grothendieck.
The operators $b$ and $B$
are all continuous and extend to the
completed topological vector spaces.
We then define $HH_\bullet(\A)$, $HC_\bullet(\A)$ and
$HP_\bullet(\A)$ for a Fr\'{e}chet algebra in exactly the same
manner as above using the complexes formed by
replacing ordinary tensor products with projective tensor
products.  Clearly the groups $HH_\bullet(\A)$, $HC_\bullet(\A)$
and $HP_\bullet(\A)$ are all functorial with respect to continuous
homomorphisms $\phi\colon \A\to \B$ of Fr\'{e}chet algebras.
It is also possible to define a (topological) \emph{bivariant}
version $HP_\bullet(\A,\B)$ of periodic cyclic theory 
\cite{CuQu,Cu1,Cu2}
however we shall make no use of this.

Finally, we summarise
some key properties of periodic cyclic homology,
thought of as a functor $HP_\bullet$ from some suitable class of
locally convex algebras to the category of graded abelian groups
(we refer to \cite{CuQu,Cu1,Cu2} for more details)

%\begin{description}
%\item[$(HP\ 1)$]
$(HP 1)$: $HP_\bullet$ is \emph{diffeotopy invariant}, i.e. if
$\alpha_0,\alpha_1 \colon \A\to \B$ are differentiably homotopic
(so that there exists a continuous homomorphism $\alpha\colon
\A\to C^\infty([0,1])\hat{\otimes}_{\pi}
\B$ with $\alpha(a)(0) = \alpha_0(a)$
and $\alpha(a)(1) = \alpha_1(a)$ for all $a\in \A$) then the
induced homomorphisms $(\alpha_0)_*, (\alpha_1)_*\colon
HP_{\bullet}(\A)\to HP_{\bullet}(\B)$ are equal.

%\item[$(HP\ 2)$]
$(HP 2)$: $HP_\bullet$ is \emph{stable}, i.e. if $\cL^1$ denotes
the
trace class operators on a separable Hilbert space $H$,
then $HP_{\bullet}(\A\hat{\otimes}_{\pi} \cL^1)$
is isomorphic to $HP_{\bullet}(\A)$.

%\item[$(HP\ 3)$]
$(HP 3)$: $HP_\bullet$ satisfies \emph{excision}, i.e. if
$0\to \A\to \B\to \mathcal{C}\to 0$ is an exact sequence of algebras
which
admits a continuous linear splitting, then we have the
six-term exact sequence
$$
\begin{array}{ccccc}
HP_0(\A) & \rightarrow & HP_0(\B)
& \rightarrow & HP_0(\mathcal{C}) \\
\uparrow \scriptstyle{\delta} & & & & \downarrow \scriptstyle{\delta} \\
HP_1(\mathcal{C}) & \leftarrow & HP_1(\B) & \leftarrow & HP_1(\A)
\end{array}
$$
%\end{description}
Of these, property $(HP 2)$ will be particularly important;
we discuss it in more detail in \S~3.3 below.

%-----------------------------------------------------------------------

\subsection{The Cyclic Homology of Smooth Functions on a Manifold}
\label{sec:one b}

We will be
interested in the Fr\'{e}chet algebra $\A = C^\infty(X)$ of smooth
functions on a manifold $X$.  Here the answer was first
computed by Connes in \cite{Connes}.

\begin{proposition}[\cite{Connes,Pfl}]
\label{prop:one}
Let $\A = C^\infty(X)$ be the Fr\'{e}chet algebra of smooth
functions on a paracompact manifold $X$.  Then
there is an isomorphism between the
(continuous) Hochschild homology groups $HH_n(\A)$ of $\A$
and the differential $n$-forms $\Omega^n(X)$
on $X$.  The cyclic homology groups $HC_n(\A)$ can
be identified with
$$
HC_n(\A) = \Omega^n(X)/\mathrm{d} \Omega^{n-1}(X)\oplus
uH^{n-2}_{\text{dR}}(X)\oplus
u^2H^{n-4}_{\text{dR}}(X)\oplus \cdots
$$ which terminates in $u^{\frac{n}{2}}H^0_{\text{dR}}(X)$ if
$n$ is even and in $u^{\frac{n-1}{2}}H^1_{\text{dR}}(X)$ if $n$ is odd.
Finally $HP_{n}
(\A) = H_{\text{dR}}^{\mathrm{ev}}(X)u^n$ if $n$ is
even and $HP_n(\A)
= H_{\text{dR}}^{\mathrm{odd}}(X)u^n$ if $n$ is odd.
\end{proposition}
In \cite{Connes} Connes proved this result
under the assumption that $X$ was compact.
Connes showed that a Kozul complex gave a topological
projective resolution of $\A$ and was able to identify the
Hochschild cohomology of $\A$ with the space of de Rham
currents on $X$.
Pflaum \cite{Pfl} later removed the
restriction that $X$ be compact.
The map $\phi_k\colon \A^{\otimes k+1} \to \Omega^k(X)$ defined by
\begin{equation}
\label{eq: CHKR map}
\phi_k(a_0\otimes a_1\otimes \cdots \otimes a_k)  = \frac{1}{k!}a_0\extd
a_1 \cdots  \extd a_k
\end{equation}
gives a map of complexes $\phi_k$ from $CC_k(\A)\otimes
\Cset((u))/u\Cset[[u]]$ to $\Omega^k(X)\otimes \Cset((u))/
u\Cset[[u]]$ where the latter complex is equipped with
the differential $ud$.  Note that $\phi_k$ extends to
the completion of $CC_k(\A)$ in the projective tensor
product topology.  Similarly, $\phi_k$ induces a map
of complexes $\phi_k\colon CC_k(\A)\otimes \Cset((u))\to \Omega^k(X)
\otimes \Cset((u))$.
It can be shown that the map $\phi_k$ induces an isomorphism
on the homology of these complexes.

%-----------------------------------------------------------------------

\subsection{Topologically H-unital algebras and topological Morita
invariance}

In this section we describe the cyclic homology
of the Schatten ideals.
We begin with the following Lemma from the thesis \cite{Baehr}
of Baehr
--- for which we are
very grateful to Joachim Cuntz for making available to us.

\begin{lemma}
Suppose $\A$ is a Fr\'{e}chet algebra such that the
multiplication map $m\colon \A\hat{\otimes}_{\pi} \A
\to \A$ has a continuous left $\A$-linear section
$s\colon \A\to \A\hat{\otimes}_{\pi} \A$.  Then
$\A$ is topologically $H$-unital.
\end{lemma}
The point is the section $s$ allows one to define a
contracting homotopy of the (topological) bar
complex of $\A$.
In particular (following \cite{Baehr}), this allows us to prove that
$\cL^1$ is topologically $H$-unital.  To see this,
recall that there is an identification $\cL^1 =
H\hat{\otimes}_{\pi} H$ under which the multiplication
map corresponds to the continuous map defined by
$a\otimes b\mapsto a(u)\otimes b^*(u)$ for a fixed
unit vector $u\in H$.  We define a continuous left
$\cL^1$-linear section of the multiplication map as follows:
\begin{equation}\label{eqn:section}
\begin{array}{rcl}
s: \cL^1 = H\hat{\otimes}_{\pi} H &\to& \cL^1 \otimes \cL^1\\[+7pt]
v\otimes w &\to& \langle \cdot,u\rangle v
\otimes \langle \cdot,w\rangle u
\end{array}
\end{equation}
where $\langle\, ,\, \rangle$
denotes the inner product on $H$.  Therefore, one concludes
that $\cL^1$ is topologically $H$-unital.  We make the remark,
following \cite{Baehr}, that $\cL^1$ is not algebraically
$H$-unital since the multiplication map restricted to the
algebraic tensor product $\cL^1\odot \cL^1$ is not surjective.
Since $\cL^1$ is topologically $H$-unital, one can use the
complex $CC_\bullet(\cL^1)$ to compute $HH_\bullet(\cL^1)$, without
first having to add a unit to $\cL^1$.  As in \cite{Lod}
we have the trace map $\tr_n\colon (\cL^1)^{\hat{\otimes}_{\pi} n}
\to \Cset^{\hat{\otimes}_{\pi} n}$ defined by
\begin{equation}
\label{eq: trace morphism}
\tr(A\otimes B \otimes \cdots \otimes C) =
\sum A_{i_0 i_1}\otimes B_{i_1 i_2} \otimes
\cdots \otimes C_{i_n i_0}
\end{equation}
where the sum is extended over all possible sets of indices.
In other words $\tr_n(A\otimes B\otimes \cdots
\otimes C) = \tr(AB\cdots C)$.
It is easy to check that $\tr_n$ defines a chain map
$\tr_\bullet\colon CC_\bullet(\cL^1)\to CC_\bullet(\Cset)$.
If $p$ is a rank one projection on $H$ with $\tr(p) = 1$
then we have a homomorphism $\mathrm{inc}\colon \Cset\to \cL^1$
which extends to define a chain map $\mathrm{inc}_\bullet
\colon CC_\bullet(\Cset)\to CC_\bullet(\cL^1)$.  Clearly
$\tr_\bullet\circ \mathrm{inc}_\bullet = \mathrm{id}$.
However it is also easy to check that the explicit chain
homotopy $h$ from $\mathrm{inc}_\bullet \circ \tr_\bullet$
to $\mathrm{id}$ given in \cite{Lod} pages 17--18
goes through in this situation.  Therefore we conclude that
$\tr_\bullet \colon HH_\bullet(\cL^1)\to HH_\bullet(\Cset)$ is
an isomorphism and hence we have the following
description of the Hochschild, cyclic and periodic cyclic
homology groups of $\cL^1$.

\begin{proposition}[\cite{Baehr,Cu3}]
\label{prop: Hoch + cyc + p cyc homology of Schatten ideals}
The Hochschild homology groups $HH_n(\cL^1)$ are all zero,
except for $n=0$ when we have $HH_0(\cL^1) = \Cset$.  It follows
therefore that
we have
$$
HC_n(\cL^1) = \begin{cases}
\Cset & \text{if}\ n\ \text{is even} \\
0 & \text{if}\ n \ \text{is odd}.
\end{cases}
$$
Finally the periodic cyclic homology groups $HP_n(\cL^1)$
equal $\Cset$ if $n$ is even and are zero if $n$ is odd.
\end{proposition}

%We remark that Baehr computes the Hochschild, cyclic
%and periodic cyclic homology of a topological algebra
%$\A$ in terms of its coinvariants \cite{Baehr}.
More generally, one can show the following
\begin{lemma}\label{thm:H-unital1}
If $\A$ is a unital
Fr\'{e}chet algebra then $\A\hat{\otimes}_{\pi} \cL^1$
is topologically $H$-unital.
\end{lemma}

The continuous
$\A\hat{\otimes}_{\pi} \cL^1$-linear section is obtained
from \eqref{eqn:section} by extending it $\A$-linearly.
In particular this Lemma has as a consequence that the
$C^\infty(X,\cL^1) = C^\infty(X)\otimes_{\pi}\cL^1$ is
topologically $H$-unital.

We can define a generalised trace
morphism $\tr_\bullet\colon CC_\bullet(\A\hat{\otimes}_{\pi}
\cL^1) \to CC_\bullet(\A)$ using a generalisation of
the formula~\eqref{eq: trace morphism} above ---
this is because we can identify $(\A\hat{\otimes}_{\pi}
\cL^1)^{\hat{\otimes}_{\pi} k} = \A^{\hat{\otimes}_{\pi} k}
\hat{\otimes}_{\pi} (\cL^1)^{\hat{\otimes}_{\pi} k}$ and
then use ~\eqref{eq: trace morphism} to define a map
$\A^{\hat{\otimes}_{\pi} k}\hat{\otimes}_{\pi}
(\cL^1)^{\hat{\otimes}_{\pi} k} \\ \to\A^{\hat{\otimes}_{\pi} k}$
which is compatible with $b$ and $B$.  The proof
of \cite{Lod} referred to above also extends to this
situation to show that $\tr_\bullet\circ \mathrm{inc}_\bullet
= \mathrm{id}$  and $\mathrm{inc}_\bullet\circ \tr_\bullet
= \mathrm{id}$ where $\mathrm{inc}_\bullet \colon
CC_\bullet(\A)\to CC_\bullet(\A\hat{\otimes}_{\pi} \cL^1)$
is induced by the homomorphism $a\mapsto a\otimes p$ where
$p$ is a projection of rank one in $H$, normalised to have
trace $\tr(p) = 1$.  Thus we have the following
generalisation of Proposition \ref{prop: Hoch + cyc + p cyc
homology of Schatten ideals}.

\begin{proposition}
\label{prop: Hoch + cyc + p cyc homology of Schatten ideals otimes A}
The generalised trace morphism $\tr_\bullet \colon
CC_\bullet(\A\hat{\otimes}_{\pi} \cL^1) \to
CC_\bullet(\A)$ induces isomorphisms
$HH_n(\A\hat{\otimes}_{\pi} \cL^1) \cong
HH_n(\A)$, $HC_n(A\hat{\otimes}_{\pi}
\cL^1)\cong HC_n(\A)$ and $HP_n(\A
\hat{\otimes}_{\pi} \cL^1)\cong HP_n(\A)$ in Hochschild,
cyclic and periodic cyclic homology respectively.
\end{proposition}

\begin{remark}
\label{rem: Hoch hom of C(X,L^1) = forms}
This Proposition can be viewed as a version of topological
Morita invariance of Hochschild, cyclic and periodic cyclic
homology, where the finite dimensional matrices are
replaced by their completion in the $\cL^1$ norm.
It also
has as a consequence that the
trace morphism $\tr_n\colon HH_n(C^\infty(X,\cL^1))
\to HH_n(C^\infty(X))$, when combined with Connes' HKR map~\eqref{eq:
CHKR map}
above, gives an identification of the Hochschild
homology of the algebra $C^\infty(X,\cL^1)$ with
the space of differential forms on $X$.
\end{remark}
We also mention the following result of
Cuntz examining the effect  in periodic cyclic homology
of replacing $\cL^1$ by $\cL^p$ for $p\ge 1$.

\begin{proposition}[\cite{Cu3}]
If $1\le p$, then the homomorphism $
\A \hat{\otimes}_{\pi} \cL^1 \to \A \hat{\otimes}_{\pi} \cL^p$
induced by the inclusion $\cL^1\hookrightarrow \cL^p$ induces
isomorphism $HP_\bullet(\A\hat{\otimes}_{\pi}
\cL^1) \stackrel{\cong}{\to} HP_\bullet(\A\hat{\otimes}_{\pi} \cL^q)$ in
periodic
cyclic homology.
\end{proposition}

When this is combined with
Proposition \ref{prop: Hoch + cyc + p cyc homology of Schatten ideals
otimes A},
we deduce that $HP_\bullet(\A) \stackrel{\cong}{\to}
HP_\bullet(\A\hat{\otimes}_{\pi}
\cL^p)$ for all $p\ge 1$, and therefore can be viewed as yet another
version of topological
Morita invariance periodic cyclic
homology, where the finite dimensional matrices are
replaced by the topological completion given by $\cL^p$.

{\em For the rest of the paper, we will assume that $X$ is a
compact manifold, even if not explicitly stated.}

The following generalizes Lemma \ref{thm:H-unital1}.

\begin{lemma} \label{thm:H-unital2}
Let $X$ be a compact manifold. Then the $*$-algebra $C^\infty(\cL^1(P))$
is topologically $H$-unital.
\end{lemma}

\begin{pf*}{Proof.}
Let $\{U_\alpha\}$ be an open cover by good open subsets of $X$,
and $\{\phi_\alpha\}$ be a partition of unity such that $\sum_\alpha
\phi_\alpha^2=1$.
Let $\cL^1(P)_\alpha = \cL^1(P))\big|_{U_\alpha}$ be the restriction,
and
$$
r_\alpha : C^\infty(\cL^1(P)) \to C^\infty(\cL^1(P)_\alpha)
$$
be the induced restriction homomorphism.
Let
$$
i_\alpha : C^\infty(\cL^1(P)_\alpha) \to C^\infty(\cL^1(P))
$$
be defined as $i_\alpha(g) = \phi_\alpha g$.
Note that $i_\alpha(r_\alpha(f)g) = fi_\alpha(g)$
for $f\in C^\infty(\cL^1(P))$.  By Proposition
\ref{prop: Hoch + cyc + p cyc homology of Schatten ideals otimes A}
and Remark \ref{rem: Hoch hom of C(X,L^1) = forms}, we know that
$C^\infty(\cL^1(P)_\alpha) = C^\infty(U_\alpha,\cL^1)$
is topologically H-unital, therefore there
is a
continuous left $C^\infty(\cL^1(P)_\alpha)$-linear section of the
multiplication map,
$$
s_\alpha : C^\infty(\cL^1(P)_\alpha) \to C^\infty(\cL^1(P)_\alpha)
\otimes_{\pi} C^\infty(\cL^1(P)_\alpha).
$$
Define the map
$$
t_\alpha : C^\infty(\cL^1(P)) \to C^\infty(\cL^1(P)) \otimes_{\pi}
C^\infty(\cL^1(P))
$$
as the composition $t_\alpha = (i_\alpha \otimes i_\alpha)\circ
s_\alpha\circ r_\alpha$.  Note that $t_\alpha$ is not
quite a section for the multiplication map, instead we
have $m\circ t_\alpha = \phi^2_\alpha t_\alpha$.
Now define $t =
\sum_\alpha t_\alpha$,
so that we get a map
$$
t : C^\infty(\cL^1(P)) \to C^\infty(\cL^1(P)) \otimes_{\pi}
C^\infty(\cL^1(P)).
$$
It is easy to check that $t$ is a continuous
left $C^\infty(\cL^1(P))$ section of the
multiplication map,
establishing the topological H-unitality of $C^\infty(\cL^1(P))$.
\end{pf*}

%-----------------------------------------------------------------------

\section{Twisted $K$-theory and twisted cohomology}
\label{sec:three}
\subsection{Twisted $K$-theory}
\label{sec:three a}

Let $X$ be a compact space which comes equipped with a
principal $PU$ bundle $P\to X$ where $PU = PU(H)$ is the
projective unitary group of a separable Hilbert space $H$.
We will assume that the topology on $PU$ is induced from the
norm topology on the full unitary group $U = U(H)$
(this has the advantage that $PU$ acquires the structure of a
Banach Lie group, this would not be true for the topology on
$PU$ induced from the strong operator topology on $U$ ---
see \S~4.4 for more details).
$PU$ acts by conjugation on the space $\K = \K(H)$ of compact
operators
on $H$.  Hence we can form the associated bundle $\K(P)$
and consider the (non-unital)
$C^*$-algebra $A = C(\K(P))$ of continuous sections of
$P \times_{PU} \K$.  The $K$-theory $K_\bullet(A)$ of the algebra
$A$ is called the \emph{twisted $K$-theory} of the pair
$(X,P)$ \cite{Ros} and is denoted $K^i(X, P)$.
The groups $K^i(X;P)$ are a functor from the category of
spaces equipped with a principal $PU$ bundle to the category
of abelian groups.  In this respect the groups $K^i(X;P)$
have many similarities with the groups $K^i(X)$, for example
there are long exact sequences associated to pairs $(X,Y)$.
In \cite{Ros} Rosenberg showed that the twisted $K$-groups
$K^i(X;P)$ could also be interpreted as vertical homotopy
classes of sections of certain bundles of Fredholm operators
associated to $P$. There are
alternative descriptions of the groups $K^i(X;P)$ ---
for example in \cite{BCMMS} the twisted $K$-groups are
interpreted in terms of bundles on $P$ twisted by a gerbe.
In \cite{Atiyah}, another equivalent definition of twisted $K$-theory 
is given,
where $PU$ is given the compact open topology instead
of the norm topology.

\subsection{Smooth subalgebras}
In this note however we want to study the $K$-theory of
the algebra $A$ from the point of view of non-commutative
geometry.  We must therefore find a smooth replacement for
$A$.   Assume from now on that $X$ is a smooth compact
manifold and that $P\to X$ is a Banach principal $PU$ bundle
(so that $P$ is a Banach manifold).

We begin by recalling some
generalities on smooth subalgebras of
$C^*$-algebras. Let $A$ be a $C^*$-algebra and
$\widetilde {A}$ be obtained by
adjoining a unit to
${A}$. Let ${\A}$ be a
$*$-subalgebra of
${A}$ and $\widetilde {\A}$ be
obtained by
adjoining a unit to ${\A}$. Then ${\A}$ is
said to be a {\em smooth subalgebra} of ${A}$ if
the
following two conditions are satisfied:

\begin{enumerate}
\item
${\A}$ is a dense $*$-subalgebra of ${A}$;

\item ${\A}$ is stable under the holomorphic
functional calculus,
that is, for any $a\in \widetilde {\A}$ and for any
function $f$ that is holomorphic in a neighbourhood
of the
spectrum of $a$ (thought of as an element in
$\widetilde
{A}$) one has $f(a) \in \widetilde {\A}$.
\end{enumerate}

Assume that ${\A}$ is a dense
$*$-subalgebra of
${A}$ such that ${\A}$ is a
Fr\'echet algebra
with a topology that is finer than that of
${A}$. A
necessary and sufficient condition for ${\A}$ to be a
smooth subalgebra is given by the {\em spectral
invariance}
condition cf. \cite[Lemma 1.2]{Schw}:

\begin{itemize}
\item $\widetilde {\A}\cap
GL(\widetilde {A}) =
GL(\widetilde {\A})$, where
$GL(\widetilde{\A})$ and
$GL(\widetilde {A})$
denote the group of
invertibles in $\widetilde{\A}$ and
$\widetilde {A}$ respectively.
\end{itemize}

Now $\cL^p$ is a dense $*$-subalgebra of $\K$. Suppose that
$1+a \in \tilde\cL^p$ is such that $1+a \in GL(\tilde \K)$, i.e.
there is $1+b  \in GL(\tilde \K)$ such that $(1+ a)(1+b) = 1$.
Then $a + b + ab = 0$, so that $b = - a - ab \in \cL^p$, i.e
$1+b  \in GL(\tilde \cL^p)$. This shows that $\cL^p$ is
a spectral invariant subalgebra of $\K$. It follows that
$\cL^p$ is a smooth dense subalgebra of $\K$.

Let  $\cL^1(P) = P
\times_{PU} \cL^1$ be the
associated algebra bundle on $X$
(here $PU$ acts on $\cL^1$ by
conjugation). Let $\A= C^\infty(\cL^1(P))$ be the
$*$-algebra of smooth sections of $\cL^1(P)$, which
is made into a Fr\'{e}chet algebra
with the topology defined in $\S 2$.

Let $\mathfrak A = C(\cL^1(P))$ be the $B^*$-algebra of continuous
sections.
$\mathfrak A$ is then a smooth, dense subalgebra of $A=  C(\K(P))$.
This can be seen as follows.
The argument given above shows that if $1+a \in \tilde{\mathfrak A}$
is such that $1+a \in GL(\tilde A)$, then $1+a \in
GL(\tilde {\mathfrak A})$.
Next, a direct generalization of the approximation argument showing
that $C^\infty(X)$
is a smooth dense subalgebra of $C(X)$ shows that $\A$ is a
smooth dense subalgebra of $\mathfrak A$. Combining these observations,
we see that $\A$ is a smooth dense subalgebra of $A$ as desired and
it also follows that $\A$ is stable under the holomorphic
functional calculus.
Therefore by a result in the appendix in \cite{Connes}, (a more detailed
proof is given in \cite{Bost}, see also chapter 3, \cite{G-BVF}) the 
inclusion
map
$\A \subset A$ induces an isomorphism
$K_i(\A) \cong K_i(A)$.

\subsection{Connes-Chern character}
Recall how the Connes-Chern character
$$\text{ch}
\colon K_i(\A)\to HP_i(\A)$$ is constructed.
Let $\mathcal C$ be a unital $*$-algebra with unit $ 1_{\mathcal C}$,
and let
${\rm Proj}(\mathcal C)$ be its set of self-adjoint projections.
Two
projections $P, Q \in {\rm Proj}(\mathcal C)$ are said to be {\it
Murray-von Neumann equivalent} if there is an element $V\in
\mathcal C$ such that $P= V^*V$ and $Q= VV^*$.
Denote $M_n(\Cset)\otimes \mathcal C= M_n( \mathcal C )$,
where $M_n(\Cset)$ denotes the square matrices of size
$n$ over $\Cset$.
Then $M_n( \mathcal C )$ is also a
$*$-algebra. Let $M_\infty( \mathcal C ) = \lim_{n\to\infty}
M_n( \mathcal C )$ be the direct limit of the embeddings of $M_n(
\mathcal C )$ in $M_{n+1}( \mathcal C )$ given by $T\to
\left(\begin{array}{cc}T & 0\\0 & 0\end{array}\right)$. Let
$V(\mathcal C) = {\rm Proj}(M_\infty( \mathcal C ))/\sim\;$ denote
the Murray-von Neumann equivalence classes of projections in
$M_\infty(\mathcal C )$.
Then $V(\mathcal C) $ is
an Abelian  semi-group under with the operation induced by the direct
sum,
and the associated Abelian group is called the Grothendieck group
$K_0(\mathcal C)$.
The homomorphism $\pi:\mathbb C \to \mathcal C$ given by $\lambda
\mapsto \lambda\cdot 1_{\mathcal C}$ induces a homomorphism
$\pi_*:K_0(\mathbb C )\cong {\mathbb Z} \to K_0(\mathcal C)$. Then
the reduced $K$-group $\tilde{K}_0(\mathcal C) $ is defined as
$\tilde{K}_0(\mathcal C)= \operatorname{coker} \pi_*\cong
K_0(\mathcal C)/{\mathbb Z}$.

Now $\mathcal A$ is a non-unital $*$-algebra; let
$\widetilde{\mathcal A}$ be the
unitization of $\A$.
By definition, the $K$-group $K_i(\mathcal A)$ is the reduced
$K$-group $\tilde{K}_i(\widetilde{\mathcal A})$ of
$\widetilde{\mathcal A}$. Explicitly,  let $P, Q$ be projections
in $M_n(\widetilde {\mathcal A} )$ such that $P-Q \in M_n(\mathcal A )$.
Then $[P-Q] \in K_0(\mathcal A)$.
Recall that the Connes-Chern character is defined as
\begin{equation}
\label{eq:11}
\mathrm{ch}([P-Q]) = \tr(P-Q) + \sum_{n\in \mathbb N}(-u)^n
\frac{(2n)!}{n!} \;
\tr\left( \left(
P-Q -\frac12 \right)
(d(P-Q))^{2n}\right),
\end{equation}
which is an element of the complex $CC_k
(\A)\otimes \Cset((u))$ with boundary
operator $b+uB$.
Now let $U$ be an invertible element in $M_n(\widetilde \A)$ such that
$U - 1 \in
             M_n(\A)$. Then $[U] \in K_1(\A)$, and the odd Connes-Chern
character is
defined as
\begin{equation}
\label{eq:12}
\mathrm{ch}([U]) =  \sum_{n\ge 0} u^n n! \;  \tr\left( \left(U^{-1} - 1
\right)
d(U -1) (d(U^{-1} -1) d(U -1) )^{n}\right),
\end{equation}
which is an element of the complex $CC_\bullet(\A)
\otimes \Cset((u))$ with boundary
operator $b+uB$.

%-----------------------------------------------------------------------

\subsection{The Dixmier-Douady Class}
\label{sec:three b}

We consider
principal bundles $P$ on $X$ with structure group $PU$.
Isomorphism classes of
such bundles on $X$ are in bijective correspondence with $H^3(X;\Zset)$.
To see this recall that by Kuiper's theorem $PU$ has the homotopy
type of a $K(\Zset,2)$ and hence $BPU$ has the homotopy type of a
$K(\Zset,3)$.  Recall that we are considering $PU$ with the
topology induced from the norm topology on $U$.
Recall that $U = U(H)$ is a Banach Lie group (see for example
\cite{Qui1}).
$U(1)$ is a closed subgroup of $U$ and hence $PU$ has a natural
structure as a Banach Lie group.  This is not the case for the
strong operator topology on $PU$; $PU$ equipped with the
strong operator topology is an example of a \emph{Polish} group.
Note that it is possible to realise $BPU$ as a
Banach manifold --- see \cite{CaCrMu} for more details.
We associate a \v{C}ech class $\delta(P)\in H^3(X;\Zset)$
to a principal $PU$ bundle $P$ on $X$ by first choosing
a good open cover $\{U_i\}_{i\in I}$ of $X$ relative to
which $P$ has $PU$-valued transition functions $g_{ij}$
(recall that a good cover is one for which every non-empty
finite intersection $U_{i_1}\cap \cdots \cap U_{i_p}$ is
contractible).
Choose lifts $\hat{g}_{ij}\colon
U_{ij}\to U$ of $g_{ij}$ and define $\epsilon_{ijk}\colon
U_{ijk}\to U(1)$ by $\hat{g}_{ij}\hat{g}_{jk}
= \hat{g}_{ik}\epsilon_{ijk}$.  Since $U(1)$ is
central in $U$ one can show that $\epsilon_{ijk}$ satisfies
the \v{C}ech $2$-cocycle condition.  Finally we choose
maps $w_{ijk}\colon U_{ijk}\to \Rset$ such that $\exp(2\pi i w_{ijk})
= \epsilon_{ijk}$ and define $n_{ijkl} = w_{jkl} - w_{ikl}
+ w_{ijl} - w_{ijk}$.  Then $n_{ijkl}$ is a \v{C}ech representative
for a class $\d(P) \in H^3(X;\Zset)$.

Freed in \cite{Freed} gives a very nice construction of a
differential $3$-form representing the image in
de Rham cohomology of the Dixmier-Douady class $\d(P)$
(note that a similar construction appeared in
\cite{Gomi} and in the context of loop groups in
\cite{MurSte}).  We recall his result here, phrased
in the language of Michael
Murray's
\emph{bundle gerbes} \cite{Mur}.
Associated to $P$ is the \emph{lifting bundle gerbe}
$L\to P^{[2]}$, where $P^{[2]} = P\times_{\pi}P$ is the
fibre product.  Recall that $L$ is the complex line associated to
the principal $U(1)$ bundle on $P^{[2]}$ obtained by pulling
back the universal $U(1)$ bundle on $PU$ via the canonical
map $g\colon P^{[2]}\to PU$.  Here $g(p_1,p_2)$ for
$(p_1,p_2)\in P^{[2]}$ is the unique element of $PU$
such that $p_2 = p_1g(p_1,p_2)$.  The product on the group $U$
induces a \emph{bundle gerbe product} on $L$, i.e. a line
bundle isomorphism which on fibres takes the form
$L_{(p_1,p_2)}\otimes L_{(p_2,p_3)}\to L_{(p_1,p_3)}$ for
points $p_1$, $p_2$ and $p_3$ of $P$ all lying in the same
fibre of $P$ over $X$.  A \emph{bundle gerbe connection} on $L$
is a connection $\nabla_L$ on the line bundle $L$ which is
preserved by the bundle gerbe product.  Thus if $s$ and
$t$ are sections of $L$ then we have $\nabla_L(st)(p_1,p_3)
= \nabla_L(s)(p_1,p_2)t(p_2,p_3) + s(p_1,p_2)\nabla_L(t)(p_2,p_3)$.
It can be shown that the curvature $F_{\nabla_L}$ of a bundle
gerbe connection $\nabla_L$ can always be written in the form
$F_{\nabla_L} = \pi_1^*f - \pi_2^*f$ for some $2$-form $f$ on $P$,
where $\pi_1$ and $\pi_2$ denote the projections onto the
first and second factors in $P^{[2]}$ respectively.  It can
also be shown that $\extd f$ is a basic form, i.e. $\extd f  =
\pi^*\omega$
for some necessarily closed $3$-form $\omega$ on $X$.  For more
details we refer to \cite{Mur}.  $\omega$ is a representative for
the image in de Rham cohomology of the Dixmier-Douady class
$\delta(P)$ of $P$.

Write $\mathfrak{u}$ and $\mathfrak{pu}$
for the Lie algebras of $U$ and $PU$ respectively.
Denote by $\text{Split}(\mathfrak{u},\mathfrak{pu})$
the set of splittings $\sigma\colon \mathfrak{pu}
\to \mathfrak{u}$ of the central extension of
Lie algebras $i\Rset \to \mathfrak{u}\stackrel{p}{\to} \mathfrak{pu}$.
Note that $\text{Split}(\mathfrak{u},\mathfrak{pu})$ is
an affine space.  $PU$ acts on $\text{Split}(\mathfrak{u},
\mathfrak{pu})$ by $g\cdot \sigma = Ad(\hat{g}^{-1})\sigma
Ad(g)$ where $\hat{g}\in U$ is such that $p(\hat{g}) = g$.
Therefore we can form the associated bundle $\text{Split}(P)$ on
$X$.  Since the fibre of $\text{Split}(P)$ is affine
we can choose a section which we denote by $\sigma$.
Observe that a splitting $s$ of the projection $\mathfrak{u}
\stackrel{p}{\to} \mathfrak{pu}$ determines a left invariant
connection $\hat{\theta}_L - s(\theta_L)$ on the principal $U(1)$ 
bundle $U\to PU$
with curvature $\frac{1}{2}[s(\theta_L),s(\theta_L)] -
\frac{1}{2}s[\theta_L,\theta_L]$ (here $\hat{\theta}_L$ and
$\theta_L$ denote the left Maurer-Cartan forms on $U$ and
$PU$ respectively).
Clearly then $\sigma$ induces a connection on the pullback
principal $U(1)$ bundle on $P^{[2]}$ and hence a connection $\nabla_L$
on the associated line bundle $L$.  $\nabla_L$ is given locally
by $\nabla_L = \extd + A_{\a}$ where $A_{\a}$ is the local
$1$-form $\hat{g}_{\a}^*\hat{\theta}_L - \pi_2^*\sigma(g^*\theta_L)$
($\hat{g}_{\a}$ is a local lift of $g\colon P^{[2]}\to PU$ to $U$).
The equivariance property
of $\sigma$ under the action of $PU$ shows that $\nabla_L$
is a bundle gerbe connection.  The curvature
$F_{\nabla_L}$ of $\nabla_L$ is easily computed
to be $F_{\nabla_L} = \extd\hat{g}^*\theta - \pi_2^*\extd
\sigma(g^*\theta)$.
Now suppose we have chosen a
connection $\Theta$ on the principal $PU$ bundle $P\to X$
with curvature $\Omega$.
A calculation, using the identity $g^*\Theta = \mathrm{Ad}(g^{-1})\Theta
+ g^*\theta_L$, shows that the $2$-form $f$ on $P$ defined
by
\begin{equation}
\label{eq:16}
f = \extd \sigma(\Theta) + \frac{1}{2}[\sigma(\Theta),
\sigma(\Theta)] - \sigma(\extd \Theta + \frac{1}{2}[\Theta,
\Theta])
\end{equation}
satisfies $F_{\nabla_L} = \pi_1^*f - \pi_2^*f$.
The results of Murray \cite{Mur} show that $\extd f = \pi^*c(P)$
for some necessarily closed $3$-form $c(P)$ on $X$.
In fact $c(P)$ is the push-forward of the basic $3$-form
on $P$ given by
\begin{equation}
\label{eq:17}
[\sigma(\Omega),\sigma(\Theta)] - \sigma([\Omega,\Theta])
- \sigma'(\Omega)
\end{equation}
where we think of $\sigma$ as a $PU$ equivariant map $P\to
\text{Split}(\mathfrak{u},\mathfrak{pu})$ and denote by $\sigma'$ its
derivative.
$c = c(P)$ is a de Rham cocycle representative for the image of
$\d(P)$ in real cohomology.
For later use we need to know how the $3$-form $c$ depends
on the choice of connection $\Theta$ on $P$ (of course there is also
a dependence on the choice of section $\sigma$ but this will not
be so important in the sequel).  When we want to avoid confusion
about the choice of connection $\Theta$ on $P$ used to construct $c$ we
will write $c = c(\Theta)$.
Any two connections $\Theta$ and $\Theta'$ on $P$ differ
by a $1$-form $A$ on $X$ with values in the bundle
$\text{ad}\, P$.  Let $f'$ be the $2$-form \eqref{eq:16}
constructed from the connection $\Theta'$ (same choice of
$\sigma$).  We have $f' = f + \pi^*\b$ where $\b$ is the
$2$-form on $X$ given by
\begin{equation}
\label{eq:18}
\b = [\sigma(\Theta),\sigma(A)] - \sigma([\Theta,A]) +
\frac{1}{2}[\sigma(A),\sigma(A)] - \frac{1}{2}\sigma([A,A])
+ \sigma'(A)
\end{equation}
Note that $\sigma'(A)$ is a scalar
valued $2$-form and therefore \eqref{eq:18} does
define a $2$-form on $X$.
It follows therefore that we have $c(\Theta') = c(\Theta) +
\extd \b$ on $X$.

%-----------------------------------------------------------------------

\subsection{Twisted Cohomology}
\label{sec:three c}

In \cite{BCMMS} (see also \cite{FHT})
it was argued that the Chern character
in twisted $K$-theory took values in a certain `twisted' cohomology
group.
If $X$ is a smooth manifold then this group can be defined as follows.
Let $c = c(P)$ be a de Rham cocycle representative for
$\d(P)$.  Thus $c$ is a closed $3$-form on $X$ with integer
periods.
Let $\Omega^\bullet(X) = \bigoplus_{n=0} \Omega^n(X)$
denote the graded algebra of differential forms on $X$.
As above we introduce a formal variable $u$ of degree $-2$
and consider the graded algebras $\Omega^\bullet(X)\otimes
\Cset((u))/u\Cset[[u]]$ and $\Omega^\bullet(X)\otimes \Cset((u))$
consisting of differential forms on $X$ with values in formal
power series in the variable $u^{-1}$ and formal Laurent
series in $u$ respectively.  We can equip both of these
graded algebras with a differential of degree one given by $\extd
- uc$.  Thus this differential is a `twisting' of
the ordinary de Rham differential $\extd$ by exterior
multiplication with the closed $3$-form $c$, scaled by the
variable $u$ so as to have degree one.
\begin{definition}
\label{def:one}
We call the homology
of the complex
$(\Omega^\bullet(X)\otimes \Cset((u)), \extd - uc)$
the \emph{twisted de Rham cohomology} of $X$ (with
respect to $c$) and denote it by $H^\bullet(X;c)$.
\end{definition}
If $c$ changes by a coboundary, $c' = c + \extd \b$, then the
complexes $(\Omega^\bullet(X)\otimes \Cset((u))/u\Cset[[u]],
\extd - uc)$ and $(\Omega^\bullet(X)\otimes
\Cset((u))/u\Cset[[u]],\extd - uc')$ are chain isomorphic
through the chain map given by multiplication with
$e^{u\b}$.  Similarly, the complexes $(\Omega^\bullet(X)
\otimes \Cset((u)),\extd - uc)$ and
$(\Omega^\bullet(X)\otimes \Cset((u)),\extd -
uc')$ are also chain isomorphic through multiplication
by $e^{u\b}$.

If $\{U_i\}_{i\in I}$ is a good cover of the manifold $X$ then
$H^\bullet(X;c)$ can be computed from the following `twisted
\v{C}ech-de Rham'
double complex

\begin{equation}
\label{eq:20}
\begin{array}{cccccc}
\!\!\!\!\!\!\!\!\!\vdots & &
\!\!\!\!\!\!\!\!\!\vdots & &
\!\!\!\!\!\!\!\!\!\vdots &            \\
\uparrow \scriptstyle{\extd - uc} & &
\uparrow \scriptstyle{\extd -uc} & &
\uparrow \scriptstyle{\extd -uc} &     \\
\displaystyle{\bigoplus_i} \Omega^1(U_i)((u)) &
\stackrel{\d}{\rightarrow} &
\displaystyle{\bigoplus_{i<j}} \Omega^1(U_{ij})((u)) &
\stackrel{\d}{\rightarrow} &
\displaystyle{\bigoplus_{i<j<k}} \Omega^1(U_{ijk})((u)) &
\stackrel{\d}{\rightarrow} \cdots            \\
\uparrow \scriptstyle{\extd - uc} & &
\uparrow \scriptstyle{\extd -uc} & &
\uparrow \scriptstyle{\extd -uc} &     \\
\displaystyle{\bigoplus_i} \Omega^0(U_i)((u)) &
\stackrel{\d}{\rightarrow} &
\displaystyle{\bigoplus_{i<j}} \Omega^0(U_{ij})((u)) &
\stackrel{\d}{\rightarrow} &
\displaystyle{\bigoplus_{i<j<k}} \Omega^0(U_{ijk})((u)) &
\stackrel{\d}{\rightarrow} \cdots            \\
\uparrow \scriptstyle{\extd - uc} & &
\uparrow \scriptstyle{\extd -uc} & &
\uparrow \scriptstyle{\extd -uc} &     \\
\displaystyle{\bigoplus_i} \Omega^{-1}(U_i)((u)) &
\stackrel{\d}{\rightarrow} &
\displaystyle{\bigoplus_{i<j}} \Omega^{-1}(U_{ij})((u)) &
\stackrel{\d}{\rightarrow} &
\displaystyle{\bigoplus_{i<j<k}} \Omega^{-1}(U_{ijk})((u)) &
\stackrel{\d}{\rightarrow} \cdots            \\
\uparrow \scriptstyle{\extd - uc} & &
\uparrow \scriptstyle{\extd - uc} & &
\uparrow \scriptstyle{\extd -uc} &     \\
\!\!\!\!\!\!\!\!\!\vdots & &
\!\!\!\!\!\!\!\!\!\vdots & &
\!\!\!\!\!\!\!\!\!\vdots &            \\
\end{array}
\end{equation}
Here by $\Omega^n(X)((u))$ we mean the elements
of $\Omega^\bullet(X)\otimes \Cset((u))$ of total
degree $n$, so that $\Omega^n(X)((u)) = \oplus_{p+q=n}
\Omega^p(X)u^q$.
Associated to this double complex are two spectral sequences
(corresponding to the two filtrations of the total complex by
rows and by columns) which compute the cohomology of the total
complex.  If we filter by rows then the $E_1$ term of the
associated spectral sequence has $E_1^{0,q}= \Omega^q(X)((u))$
and $E_1^{p,q} = 0$ for $p>0$.
Therefore the spectral
sequence collapses at the $E_2$ term and computes the twisted
de Rham cohomology of $X$.  However if we filter by columns
then the resulting spectral sequence has $E_2^{p,q}= H^p(X)u^q$ for
$q$ even and $E_2^{p,q} = 0$ for $q$ odd.  Further one can show
that the differential $d_3\colon E_3^{p,q}\to E_3^{p+3,q-2}$ is
cup product with $[c]u$.  We summarise this discussion in the
following proposition.
\begin{proposition}
\label{prop:three}
There exists a spectral sequence converging to
$H^\bullet(X;c)$ with $E_2^{p,q} =
H^p(X;\Cset)u^q$ if $q$ is even and $E_2^{p,q} = 0$
if $q$ is odd.  The even differentials $d_{2r}$ are
all zero and the differential $d_3\colon
E_3^{p,q}\to E_3^{p+3,q-2}$ is cup product with
$[c]u$.
\end{proposition}
Note that a similar spectral sequence was obtained in
\cite{FHT} by different means.
\begin{remark}
In this next section, we are going to relate periodic
cyclic homology to twisted de Rham cohomology.  For
this reason we like to think of twisted de
Rham cohomology as the \emph{homology} groups
of a complex and we will therefore equip the complex $\Omega^\bullet(X)
\otimes \Cset((u))$ with a differential of degree $-1$
given by $u\extd - u^2c$.  The only effect
of this change is that the homology groups
$H_p(\Omega^\bullet(X)\otimes \Cset((u)))$ of
the complex $\Omega^\bullet(X)\otimes \Cset((u))$ equipped
with the differential $u\extd -
u^2c)$ are equal to $u$ times the
twisted cohomology groups $H^p(X,c)$.
For this reason we will still refer to
$H_p(\Omega^\bullet(X)\otimes \Cset((u)))$ as the
twisted de Rham cohomology of $X$.
\end{remark}

%We conjecture that all
%further differentials after $d_3$ are zero.
%\begin{remark}
%Recall that in \cite{Atiyah} Atiyah and Segal defined
%a twisted cohomology theory by noting that cup product with
%$[c]$ defined a differential on the graded algebra $H(X)$.
%They therefore obtain a graded cohomology theory $H(X;[c])
%= \oplus_{i=0}H^i(X;[c])$ with
%$$
%H^i(X;[c]) = \frac{\ker\{ \cup[c] \colon H^i(X)\to H^{i+3}(X)\}}
%{\text{im}\{ \cup[c]\colon H^{i-3}(X) \to H^i(X)\}}
%$$
%where we set $H^i(X) = 0$ for $i<0$.  The twisted cohomology
%theory of Atiyah and Segal has the advantage over the twisted
%de Rham theory described above in that it is independent of
%the choice of cocycle representative for $\d(P)$.
%Note however that if our conjecture above is true then the
%associated graded group $GH(X;c)$ would be isomorphic
%to the twisted cohomology of Atiyah and Segal
%and hence $H(X;c)$ would be isomorphic (as a vector space)
%to $H(X;[c])$.  We are indebted to Prof. G. Segal
%for kindly providing us with the following
%demonstration that this is indeed the case.
%Define a filtration on the complex $\Omega^\bullet(X)$
%by defining subcomplexes $F_p$ by $F_p =
%\oplus_{q\geq p}\Omega^q(X)$.  The spectral sequence associated
%to this filtration has $E_2$ term equal to the cohomology
%of $H(X)$ of $X$.  The differential $d_2$ is zero, the
%differential $d_3$ is cup product with $[c]$ and all
%further differentials are zero.  This implies the result above.
%\end{remark}

%-----------------------------------------------------------------------

\section{Periodic Cyclic Homology and Twisted Cohomology}
\label{sec:four}
\subsection{A Chain Map}
\label{sec:four a}

In \cite{Gor1,Gor2} Gorokhovsky introduces the notion of a
generalised cycle and defines the character of such an object,
which turns out to be a cocycle in the $(b,B)$-bicomplex
computing cyclic cohomology.  We will use Gorokhovsky's
formalism to construct a map from the $(b,B)$-bicomplex
computing cyclic homology to the twisted de Rham
complex.

Define bundles $\cL^1(P)$ and $\B(P)$ associated to $P$ via the
adjoint action of $PU$ on $\cL^1$ and $\B(H)$ respectively.
          From now on let us denote by $\A$ the Fr\'{e}chet algebra
$C^\infty(\cL^1(P))$.
Choose a connection $\Theta$ on the principal bundle $P$
and let $c = c(\Theta)$ be the associated
$3$-form~(\ref{eq:17}).
Define a connection $\nabla$ on
$\cL^1(P)$ by $\nabla(s) = \extd s + [\sigma(\Theta),s]$ for a section
$s\in C^{\infty}(\cL^1(P))$.  $\nabla$ acts as a
derivation of degree $1$ on the graded algebra $\Omega^\bullet(\cL^1(P))
= \bigoplus_{n=0}\Omega^n(\cL^1(P))$
of forms on $X$ with values in $\cL^1(P)$.  Let $\Omega$ be the
curvature of the connection $\Theta$ on the principal $PU$
bundle $P\to X$.  $\sigma(\Omega)$ can be identified as a section
of the bundle $\Omega^2(\B(P)))$.  Since $\cL^1$ is an ideal in
$\B(H)$, there is an action of $\Omega^\bullet(\B(P))$ on
$\Omega^\bullet(\cL^1(P))$.  In particular, $\sigma(\Omega)$
acts as a multiplier on the graded algebra
$\Omega^\bullet(\cL^1(P))$.
Note that $\nabla(\sigma(\Omega)s) =
\sigma(\Omega)\nabla(s) -cs$ where $c$ acts on $s$ by the
obvious action of $\Omega^\bullet(X)$ on $\Omega^\bullet(\cL^1(P))$.
Here $c$ is the $3$-form representing the Dixmier-Douady class of
the bundle $P$ obtained in~\eqref{eq:17}.
Finally observe that the ordinary operator trace $\tr$
defines a trace map $\tr\colon \Omega^\bullet(\cL^1(P)) \to
\Omega^\bullet(X)$
satisfying the properties that $\tr(s_1s_2) =
(-1)^{|s_1||s_2|}\tr(s_2s_1)$ for $s_1,s_2\in
\Omega^\bullet(\cL^1(P))$ and $\tr(\nabla(s)) = \extd
\tr(s)$ for all $s\in \Omega^\bullet(\cL^1(P))$.

\begin{definition}
\label{def:two}
Following Gorokhovsky \cite{Gor1} define a map
$\mathrm{Ch}_k \colon CC_k(\A)\to \Omega^\bullet(X) \otimes
\Cset((u))$ by the JLO-type formula
\begin{equation}
\label{eq: gen CHKR map}
\Ch(\tilde{a}_0,a_1,\ldots,a_k) =
\int_{\Delta^k}\tr(\tilde{a}_0e^{-s_0\sigma(\Omega)u}
\nabla(a_1)\cdots \nabla(a_k)e^{-s_k\sigma(\Omega)u})
ds_1\cdots ds_k.
\end{equation}
Note that $\Ch_k$ is of degree zero.
In addition we also put $\Ch_0(a_0) =
\tr(a_0e^{-\sigma(\Omega)u})$.  $\Ch_k$ extends to
a $\Cset((u))$-module map $\Ch\colon CC_\bullet(\A)\otimes
\Cset((u))\to \Omega^\bullet(X)\otimes \Cset((u))$ and descends to the
quotient to give a map (also written $\Ch$)
$\Ch\colon CC_\bullet(\A)\otimes \Cset((u))/u\Cset[[u]]\to
\Omega^\bullet(X)\otimes \Cset((u))/u\Cset[[u]]$.  Finally, setting
$u=0$ defines a map $\Ch\colon CC_\bullet(\A)\to \Omega^\bullet(X)$
given by
\begin{equation}
\label{eq: Hoch gen CHKR map}
\Ch(\tilde{a}_0,a_1,\ldots,a_k) = \frac{1}{k!}\tr(\tilde{a}_0
\nabla(a_1)\cdots \nabla(a_k)).
\end{equation}
\end{definition}

\begin{remark}
If the bundle $P$ is trivial, and we take for the connection
$\nabla$ on $\cL^1(P)$ the trivial connection $\extd$ then the
map $\Ch$ reduces to the expression~\eqref{eq: CHKR map}.
\end{remark}

\begin{remark}
Note that the map $\Ch\colon CC_\bullet(\A)\otimes \Cset((u))
\to \Omega^\bullet(X)\otimes \Cset((u))$ depends on a choice
of connection $\nabla$ on the bundle $\cL^1(P)$.  When
we want to make this precise we will write $\Ch =
\Ch(\nabla)$.
\end{remark}

\begin{remark}
In the formula~\eqref{eq: gen CHKR map} above we understand that the
expression appearing as the integrand,
$\tilde{a}_0e^{-s_0\sigma(\Omega)u}
\nabla(a_1) \cdots \nabla(a_k)e^{-s_k\sigma(\Omega)u}$,
is to mean
$$
a_0e^{-s_0\sigma(\Omega)u}\nabla(a_1)\cdots
\nabla(a_k)e^{-s_k\sigma(\Omega)u} + \lambda e^{-s_0
\sigma(\Omega)u}\nabla(a_1)\cdots \nabla(a_k)e^{-s_k
\sigma(\Omega)u}
$$
if $\tilde{a}_0 = (a_0,\lambda)$.
\end{remark}

Each of the maps $\Ch$ defined above is a morphism
of complexes.  To see this we need the following
Proposition.
\begin{proposition}
For each of the maps $\Ch$ defined above, we have
\begin{equation}
\label{eq: gen CHKR map is a chain map}
\Ch\circ(b+uB) = (u\mathrm{d}-u^2c)\circ \Ch.
\end{equation}
\end{proposition}

\begin{pf*}{Proof.}
We first compute $\extd \Ch(\tilde{a}_0,a_1,\ldots,
a_k)$.  This equals
\begin{multline*}
\int_{\Delta^k}\tr(\nabla(a_0)e^{-s_0\sigma(\Omega)u}
\nabla(a_1) \cdots \nabla(a_k) e^{-s_k\sigma(\Omega)u})
ds_1\cdots ds_k                                          \\
+ cu\int_{\Delta^k}\tr(a_0e^{-s_0\sigma(\Omega)u}
\nabla(a_1)\cdots \cdots \nabla(a_k)e^{-s_k\sigma(\Omega)u})
ds_1\cdots ds_k                                           \\
+ \sum^k_{i=1}(-1)^{i-1}\int_{\Delta^k}\tr
(\tilde{a}_0e^{-s_0\sigma(\Omega)u}\nabla(a_1)
\cdots e^{-s_{i-1}\sigma(\Omega)u}[\sigma(\Omega),a_i]        \\
e^{-s_i\sigma(\Omega)u}\nabla(a_{i+1})
\cdots \nabla(a_k)e^{-s_k\sigma(\Omega)u})ds_1\cdots ds_k,
\end{multline*}
where we use the fact that $\nabla(\sigma(\Omega)) = -c$
and $\nabla^2(a) = [\sigma(\Omega),a]$.  Therefore we
have that $(u\extd - u^2c)\circ \Ch(\tilde{a}_0,
a_1,\ldots,a_k)$ is equal to the sum of the
two terms,
\begin{equation}
\label{eq: first term}
u\int_{\Delta^k}\tr(\nabla(a_0)e^{-s_0\sigma(\Omega)u}
\nabla(a_1)\cdots \nabla(a_k)e^{-s_k\sigma(\Omega)u})
ds_1\cdots ds_k,
\end{equation}
and
\begin{multline}
\label{eq: second term}
\sum^k_{i=1}(-1)^{i-1} u\int_{\Delta^k}
\tr(\tilde{a}_0e^{-s_0\sigma(\Omega)u}
\nabla(a_1) \cdots e^{-s_{i-1}\sigma(\Omega)u}[\sigma(\Omega),a_i]
e^{-s_i\sigma(\Omega)u}                              \\
\nabla(a_{i+1}) \cdots
         \nabla(a_k)e^{-s_k\sigma(\Omega)u})ds_1\cdots ds_k.
\end{multline}
The first expression~\eqref{eq: first term}
is equal
to
$(\Ch\circ uB)(\tilde{a}_0,a_1,\ldots,a_k)$.  To see
this, we compute $(\Ch\circ uB)(\tilde{a}_0,
a_1,\ldots,a_k)$.  We get that this is equal to
\begin{multline*}
\int_{\Delta^{k+1}}\tr(e^{-s_0\sigma(\Omega)u}\nabla(a_0)
e^{-s_1\sigma(\Omega)u} \nabla(a_1)\cdots \nabla(a_k)e^{-s_{
k+1}\sigma(\Omega)u})ds_1\cdots ds_{k+1}           \\
+ \sum^k_{i=1}(-1)^{ik} \int_{\Delta^{k+1}}
\tr(e^{-s_0\sigma(\Omega)u}\nabla(a_i) e^{-s_1\sigma(\Omega)u}
\nabla(a_{i+1}) \cdots \nabla(a_k)e^{-s_{k-i+1}\sigma(\Omega)u} \\
\nabla(a_0)e^{-s_{k-i+2}\sigma(\Omega)u} \cdots
\nabla(a_{i-1})e^{-s_{k+1}\sigma(\Omega)u})ds_1\cdots ds_{k+1}.
\end{multline*}
Commuting differential forms inside the trace gives us that
this is equal to
\begin{multline*}
\int_{\Delta^{k+1}}\tr(\nabla(a_0)e^{-s_1\sigma(\Omega)u}
\nabla(a_1)e^{-s_2\sigma(\Omega)u} \cdots \nabla(a_k)
e^{-(s_0+s_{k+1})\sigma(\Omega)u})ds_1\cdots ds_{k+1} \\
+ \sum^k_{i=1}\int_{\Delta^{k+1}}\tr(\nabla(a_0)
e^{-s_{k-i+2}\sigma(\Omega)u} \cdots
\nabla(a_{i-1})e^{-(s_0+s_{k+1})\sigma(\Omega)u}
\nabla(a_i)e^{-s_1\sigma(\Omega)u}               \\
\nabla(a_{i+1})e^{-s_2\sigma(\Omega)u}
\cdots \nabla(a_k)e^{-s_{k-i+1}\sigma(\Omega)u})
ds_1\cdots ds_{k+1}.
\end{multline*}
This is then equal to
\begin{multline*}
\sum^k_{i=1}\int_{\Delta^k}
s_{i-1} \tr(\nabla(a_0)e^{-s_0\sigma(\Omega)u} \nabla(a_1)
e^{-s_1\sigma(\Omega)u} \cdots \nabla(a_k)
e^{-s_k\sigma(\Omega)u})ds_1\cdots ds_k                  \\
+ \int_{\Delta^k}s_k\tr(\nabla(a_0)e^{-s_0\sigma(\Omega)u}
\nabla(a_1) \cdots \nabla(a_k)e^{-s_k\sigma(\Omega)u})ds_1\cdots ds_k,
\end{multline*}
which is equal to~\eqref{eq: first term}.
We now show that the second expression~\eqref{eq: second term}
can be identified
with $(\Ch\circ b)(\tilde{a}_0,a_1,\ldots,a_k)$.
To do this, we compute
$\Ch b(\tilde{a}_0,a_1,\ldots,a_k)$.  This is equal to the
sum
\begin{multline*}
\sum^k_{i=1}(-1)^{i-1}\int_{\Delta^{k-1}}\tr(\tilde{a}_0
e^{-s_0\sigma(\Omega)u}\nabla(a_1) \cdots e^{-s_{i-2}
\sigma(\Omega)u} \nabla(a_{i-1})                \\
[a_i,e^{-s_{i-1}\sigma(\Omega)u}]\nabla(a_{i+1})
\cdots \nabla(a_k)e^{-s_{k-1}\sigma(\Omega)u})
ds_1\cdots ds_{k-1}
\end{multline*}
As in \cite{Qui} Equation~7.2 we have the
general differentiation formula
$$
D(e^K) = \int^1_0 e^{(1-s)K}D(K)e^{sK}ds
$$
where $D$ is a derivation.  Using this formula
for the derivations $[a_i,\cdot\,]$ this equals
\begin{multline*}
u\sum^k_{i=1}(-1)^{i+1}\int_{\Delta^{k-1}}\int^1_0
s_{i-1}\tr(\tilde{a}_0
e^{-s_0\sigma(\Omega)u}\nabla(a_1)\cdots e^{-s_{i-2}
\sigma(\Omega)u} \nabla(a_{i-1}) \\
e^{-s_{i-1}(1-t)\sigma(\Omega)u}
[\sigma(\Omega),a_i]e^{-s_{i-1}t\sigma(\Omega)u}
\cdots \nabla(a_k)e^{-s_{k-1}\sigma(\Omega)u})dtds_1\cdots ds_{k-1}
\end{multline*}
Changing variables in the integral shows that this expression
is equal to~\eqref{eq: second term}.
\end{pf*}
\begin{remark}
By looking at the part of equation~\eqref{eq: gen CHKR map
is a chain map} which is of degree zero in $u$, we see
that the map~\eqref{eq: Hoch gen CHKR map} is a chain map
from the Hochschild complex $CC_\bullet(\A)$ to the
following complex
\begin{equation}
\label{eq: Hoch de Rham complex}
0\to \Omega^0(X)\stackrel{0}{\to}
\Omega^1(X) \stackrel{0}{\to}
\Omega^2(X)\stackrel{0}{\to} \cdots
\end{equation}
\end{remark}

%-----------------------------------------------------------------------

\subsection{A Homotopy Formula}
\label{sec:four b}

We now want to prove a homotopy formula showing what happens
when we perturb the connection $\nabla$.  Recall that $\nabla$
depended on a choice of connection $\Theta$ on the principal
bundle $P$ and a choice of section $\sigma$ of the associated
bundle $\mathrm{Split}(P)$ of splittings.  Changing $\Theta$ to a
new connection $\Theta'$ will change the twisting cocycle
$c(\Theta)$ by a coboundary; i.e. $c' = c+ \extd \b$
where $c' = c(\Theta')$ and $c=c(\Theta)$ and hence changes the
target space of the map $\mathrm{Ch}$.  Recall though that the
complexes $(\Omega^\bullet(X)\otimes \Cset((u)),u\extd - u^2c)$
and
$(\Omega^\bullet(X)\otimes \Cset((u)),u\extd -
u^2c')$ are isomorphic through the
chain map
$e^{u\b}$.

Let $\Theta$ and $\Theta'$ be two connections on the
prinicipal $PU$ bundle $P\to X$.  Then $\Theta$ and $\Theta'$
are related by a $1$-form $A$ with values in the bundle
$\ad P$: $\Theta' = \Theta + A$.  Form the family of
connections $\Theta_t = \Theta + tA$ for $0\leq t\leq 1$.  Let
$\Omega_t$
denote the curvature of the connection $\Theta_t$.
Denote also by $c_t$ the $3$-form~(\ref{eq:17})
associated to the connection $\Theta_t$.  The family $\Theta_t$
induces a family of connections $\nabla_t$ on the bundle
$\cL^1(P)$ by the same formula as above, $\nabla_t(s) =
\extd s + [\sigma(\Theta_t),s]$ for $s$ a section of
$\cL^1(P)$.  We have $\nabla_t(\Omega_t) = c_t$.
Define a $2$-form $\b_t$ on $X$ by the
formula
\begin{equation}
\label{eq:30}
\b_t = t([\sigma(\Theta),\sigma(A)] - \sigma([\Theta,A])) +
\frac{t^2}{2}([\sigma(A),\sigma(A)] - \sigma([A,A])) +
t\sigma'(A).
\end{equation}
Note that we have $f_t = f + \pi^*\b_t$ and hence $c_t = c +\extd \b_t$.
Define a connection $1$-form $\tilde{\Theta}$ on the
principal $PU$ bundle $P\times [0,1] \to X\times [0,1]$
by $\tilde{\Theta} = \Theta_t$.  The curvature
of the connection $\tilde{\Theta}$ is
$\tilde{\Omega} = \Omega_t + \extd tA $.
If we apply the section $\sigma$
to $\tilde{\Omega}$ we get $\sigma(\tilde{\Omega}) =
\sigma(\Omega_t) + \extd t \sigma(A)$.
Denote by $\tilde{\nabla}$ the connection $\nabla_t + \extd t
\partial_t$ on the associated bundle $\mathcal{L}^1(P\times [0,1])$
(we will denote differentiation with respect to $t$ by
$\partial_t$ and sometimes by a dot).
An easy calculation shows that $\tilde{\nabla}\sigma(\tilde{\Omega}) =
c + \extd \b_t + \extd t \dot{\b}_t$.

Now we are ready for the homotopy formula.  Let us write
$\tilde{\A}$ for the algebra $C^\infty(\cL^1(P\times[0,1]))$
and $\tilde{\Ch} = \Ch(\tilde{\nabla})$.
We have a natural map $\rho\colon\A\to \tilde{\A}$ which sends
a section of $\cL^1(P)$ to a section of $\cL^1(P\times[0,1])$
constant in the $t$ direction.  $\rho$ induces a map
(also denoted $\rho$) $CC_k(\A)\to CC_k(\tilde{A})$ and
hence also maps $\rho\colon CC_k(\A)\otimes \Cset((u))\to
CC_k(\tilde{\A})\otimes \Cset((u))$ and
$\rho\colon CC_k(\A)\otimes \Cset((u))/u\Cset[[u]]\to
CC_k(\tilde{\A})\otimes \Cset((u))/u\Cset[[u]]$.
Notice that we have $\tilde{\nabla}
\rho(a) = \nabla_t\rho(a)$ for $a\in \A$ and hence
$\tilde{\Ch}
\rho|_{t=1} = \Ch(\nabla')$,
$\tilde{\Ch}
\rho|_{t=0} = \Ch(\nabla)$.
Define a map $K\colon CC_k(\A)\to \Omega^\bullet(X)\otimes
\Cset((u))$ of degree one by
\begin{equation}
\label{eq: homotopy operator}
K = \int^1_0 u^{-1}e^{u\b_t}\iota_{\partial_t}\tilde{\Ch}\, dt.
\end{equation}
Here $\iota_{\partial_t}$ denotes contraction with the
vector field $\partial_t$ on $[0,1]$.  $K$ induces maps,
also denoted $K$, $K\colon CC_k(\A)\otimes \Cset((u))\to
\Omega^\bullet(X)\otimes \Cset((u))$ and $K\colon
CC_k(\A)\otimes \Cset((u))/u\Cset[[u]]\to \Omega^\bullet(X)
\otimes \Cset((u))/u\Cset[[u]]$.
We have the following Proposition.
\begin{proposition}
\label{prop: gen CHKR map indep of choice of connection}
The maps $K\colon CC_\bullet(\A)\otimes \Cset((u))\to
\Omega^\bullet(X)\otimes \Cset((u))$ and $K\colon
CC_\bullet(\A)\otimes \Cset((u))/u\Cset[[u]]\to
\Omega^\bullet(X)\otimes \Cset((u))/u\Cset[[u]]$ defined
above are chain homotopies --- we have the formula
\begin{equation}
\label{eq: homotopy formula}
e^{u\b}\Ch(\nabla') - \Ch(\nabla) = K\circ (b+uB) -
(u\mathrm{d} - u^{2}c)\circ K.
\end{equation}
\end{proposition}
\begin{pf*}{Proof.}
       From equation~\eqref{eq: gen CHKR map is a chain map}
we have $\tilde{\Ch}\circ (b+uB)
= (u(\extd + \extd t\partial_t)
- u^{2}(c + \extd \b_t + \extd t\dot{\b}_t)
\circ \tilde{\Ch}$.  Multiply both sides of this
equation by $e^{u\b_t}$.  Then we get that
$$
e^{u\b_t}\tilde{\Ch}\circ (b+uB) =
(u\extd - u^{2}c)\circ e^{u\b_t}\tilde{\Ch}
+ u\extd t\partial_t e^{u\b_t}\tilde{\Ch}.
$$
Multiplying both sides now by $u^{-1}$ and integrating from
$t=0$ to $t=1$ gives the result.
\end{pf*}
\begin{remark}
\label{rem: homotopy invariance of Hoch gen CHKR map}
If we look at the part of equation~\eqref{eq: homotopy
formula} which is of degree zero in $u$, then we find that
we have $\Ch(\nabla') - \Ch(\nabla) = K_0\circ b$ where
$\Ch(\nabla')$ and $\Ch(\nabla)$ are the maps given
in equation~\eqref{eq: Hoch gen CHKR map}.  Here $K_0$ is
the degree zero in $u$ component of the map~\eqref{eq: homotopy
operator} above so that
$$
K_0 = -\int^1_0 \b_t \iota_{\partial_t}\tilde{\Ch}\, dt.
$$
This shows that the chain map $\Ch_\bullet\colon CC_\bullet(\A)\to
\Omega^\bullet(X)$ from the Hochschild complex to the
complex~\eqref{eq: Hoch de Rham complex}
induces a map on homology groups which is independent
of the choice of connection $\nabla$.
\end{remark}

%-----------------------------------------------------------------------

\subsection{Proof of the Main Result}
\label{sec:four c}

  From now on, unless we explicitly state otherwise,
we denote the completed projective tensor product
$\hat{\otimes}_{\pi}$ by $\otimes$.
Let $\{U_i\}_{i\in I}$ be a good open cover of the
compact manifold $X$.  Let $\A_i$ denote the algebra
$C^\infty(\cL^1(P)|_{U_i})$ of sections of $\cL^1(P)$ over $U_i$,
$\A_{ij}$ denote the algebra $C^\infty(\cL^1(P)|_{U_{ij}})$
of sections of $\cL^1(P)$ over $U_{ij}$ and so on.
We can identify the algebra $\oplus_i \A_i$ with the algebra
$C^\infty(\sqcup \cL^1(P)|_{U_i},\sqcup U_i)$ of smooth sections
of the bundle $\sqcup \cL^1(P)_{U_i}$ over the disjoint union
$\sqcup U_i$.
We can make similar identifications of $\oplus \A_{ij}$ and so on.
The sequence
\begin{equation}
\label{eq:37}
0\to \A \stackrel{\delta}{\to} \bigoplus_i \A_i
\stackrel{\delta}{\to}
\bigoplus_{i,j} \A_{ij}\stackrel{\delta}{\to} \cdots
\end{equation}
is then exact where the map $\delta$ is defined as follows.
Suppose that $\underline{s} = (s_{i_1\dots i_p})
\in \oplus_{i_1,\ldots,i_p}\A_{i_1\cdots i_p}$.  Then we put
$$
\delta(\underline{s})_{i_1\cdots i_{p+1}} =
s_{i_2\cdots i_{p+1}}|_{U_{i_1\cdots i_{p+1}}} -
s_{i_1 i_3\cdots i_{p+1}}|_{U_{i_1\cdots i_{p+1}}} +
\cdots + (-1)^{p+1}s_{i_1\cdots i_p}|_{U_{i_1\cdots i_{p+1}}}.
$$
Another way to think of $\delta$ is the following.  There
are $p+1$ inclusions $\sqcup U_{i_1\cdots i_{p+1}} \to
\sqcup U_{i_1\cdots i_p}$ which induce $p+1$ restriction
homomorphisms $d_i^*\colon \oplus_{i_1,\ldots,i_p}
\A_{i_1\cdots i_{p+1}}\to \oplus_{i_1,\dots,i_p}
\A_{i_1\cdots i_p}$.  $\delta$ is then the alternating
sum $\delta = \sum (-1)^i d_i^*$.  To see that the
sequence above is exact, suppose that $\delta(\underline{s}) = 0$.
Define $\underline{t} = (t_{i_1\cdots i_{p-1}}) \in
\oplus_{i_1,\ldots,i_{p-1}}\A_{i_1\cdots i_{p-1}}$ by
$t_{i_1\cdots i_{p-1}} = (-1)^p\sum_{i_p} \rho_{i_p} s_{i_1\cdots
i_p}$ for some partition of unity $\{\rho_i\}_{i\in I}$
subordinate to the open cover $\{U_i\}_{i\in I}$.  For
each $i_p\in I$ such that $U_{i_1\cdots i_p}\neq \emptyset$,
$(-1)^p\rho_{i_p}s_{i_1\cdots i_p}$ defines a section of
$\cL^1(P)|_{U_{i_1\cdots i_{p-1}}}$.  It is then easy to see
that $\delta(\underline{t})_{i_1\cdots i_p} =
s_{i_1\cdots i_p}$.

This exactness argument generalises for the algebra
$\A^{\otimes n}$.
More precisely the sequence
\begin{equation}
\label{eq:38}
0\to \A^{\otimes n} \stackrel{\delta}{\to}
(\bigoplus_i \A_i)^{\otimes n} \stackrel{\delta}{\to}
(\bigoplus \A_{ij})^{\otimes n} \stackrel{\delta}{\to}
\cdots
\end{equation}
is exact.  To see this, first observe that
$$
(\bigoplus_i \A_i)^{\otimes n} =
\bigoplus_{(i_1,\ldots,i_n)} \A_{i_1}\otimes
\cdots \A_{i_n}  = \bigoplus_{(i_1,\ldots,i_n)}
C^\infty(\cL^1(P)^{\boxtimes n},U_{i_1}\times \cdots U_{i_n}).
$$
We have similar identifications of the other tensor
products
$(\oplus_{i_1,\ldots,i_p}\A_{i_1\cdots i_p})^{\otimes n}$;
for example $(\oplus_{i,j}\A_{ij})^{\otimes n}$ is equal
to
$$
\bigoplus_{(i_1,\ldots,i_n),(j_1,\ldots,j_n)}
C^\infty(\cL^1(P)^{\boxtimes n},(U_{i_1}\times \cdots
\times U_{i_n})\cap (U_{j_1}\times \cdots \times
U_{j_n})).
$$
Recall also that $\A^{\otimes n} = C^\infty(\cL^1(P)^{\boxtimes n} ,
X^n)$.  Therefore we are in the situation above with
$X$ replaced by $X^n$ and the cover $\{U_i\}_{i\in I}$
of $X$ replaced by the cover $\{U_{i_1}\times \cdots \times
U_{i_n}\}_{(i_1,\ldots,i_n)\in I^n}$ of $X^n$.
The same partition of unity argument used above
shows that the sequence is exact.  Let us describe in more
detail the map $\delta$.  First we need some notation.  Let us
write $\underline{i} = (i_1,\ldots,i_n)$ and $U_{\underline{i}} =
U_{i_1}\times \cdots \times U_{i_n}$.  Then
$(\underline{i}_1,\ldots,\underline{i}_p) =
((i^1_1,\ldots,i_1^n),(i_2^1,\ldots,i_2^n),\ldots ,
(i_p^1,\ldots,i_p^n))$.  The $p+1$ injections
$\sqcup U_{i_1\cdots i_{p+1}}\to \sqcup
U_{i_1\cdots i_p}$ induce $p+1$ injections
$\sqcup U_{\underline{i}_1}\cap \cdots \cap
U_{\underline{i}_{p+1}} \to \sqcup U_{\underline{i}_1}
\cap \cdots \cap U_{\underline{i}_p}$ and hence $p+1$
restriction homomorphisms
$$
d_i^*\colon \bigoplus_{\underline{i}_1,\ldots,\underline{i}_p}
C^\infty(\cL^1(P)^{\boxtimes n},U_{\underline{i}_1}
\cap \cdots \cap U_{\underline{i}_p})\to
C^\infty(\cL^1(P)^{\boxtimes n},U_{\underline{i}_1}
\cap \cdots \cap U_{\underline{i}_{p+1}}).
$$
$\delta$ is again the alternating sum $\sum (-1)^i d_i^*$.
We have the following Proposition.
\begin{proposition}
\label{prop: exact seq of complexes}
We have exact sequences of complexes
\begin{align}
& 0\to CC_{\bullet}(\A) \stackrel{\d}{\to}
CC_{\bullet}(\bigoplus_i\A_i) \stackrel{\d}{\to}
CC_{\bullet}(\bigoplus_{i<j}\A_{ij}) \stackrel{\d}{\to}
\cdots  \label{eq:39}                                  \\
& 0 \to CC_{\bullet}(\A)\otimes
\frac{\Cset((u))}{u\Cset[[u]]} \stackrel{\d}{\to}
CC_{\bullet}(\bigoplus_i\A_i)\otimes
\frac{\Cset((u))}{u\Cset[[u]]} \stackrel{\d}{\to}
CC_{\bullet}(\bigoplus_{i<j}\A_{ij})\otimes
\frac{\Cset((u))}{u\Cset[[u]]} \stackrel{\d}{\to} \cdots
\label{eq:40}  \\
& 0 \to CC_\bullet(\A)\otimes\Cset((u)) \stackrel{\d}{\to}
CC_\bullet(\bigoplus_i\A_i)\otimes\Cset((u)) \stackrel{\d}{\to}
CC_\bullet(\bigoplus_{i<j}\A_{ij})\otimes\Cset((u))
\stackrel{\d}{\to} \cdots \label{eq:41}
\end{align}
\end{proposition}

\begin{pf*}{Proof.}
To prove this we  have to show that
the map $\d\colon (\oplus
\A_{i_1\cdots i_p})^{\otimes n}\to
(\oplus\A_{i_1
\cdots i_{p+1}})^{\otimes n}$
is compatible with the maps $b'$, $b$, $\lambda$ and $B$
of Section~\ref{sec:one} .  But $b'$, $b$, $\lambda$ and $B$
commute with all the homomorphisms $d_i^*$ and
hence commute with $\delta$.
\end{pf*}
This Proposition allows us to give a `\v{C}ech-de Rham
style' proof of the following theorem.
\begin{theorem}
\label{thm: main result}
Let $X$ be a compact manifold and let $P\to X$ be a
principal $PU$ bundle with connection $1$-form
$\Theta$.  Let $c = c(\Theta)$ denote the de Rham
representative of the characteristic
class in $H^3(X;\Zset)$ associated to $P$.  Then the
following results are true:
\begin{enumerate}
\item The chain map $\mathrm{Ch}$ from the Hochschild
complex $CC_\bullet(\A)$ to the
complex~\eqref{eq: Hoch de Rham complex}
is a
quasi-isomorphism.

\item The chain map $\mathrm{Ch}$ from the cyclic
complex $CC_\bullet(\A)\otimes \Cset((u))/u\Cset[[u]]$
for $\A$
to the complex $(\Omega^\bullet(X)\otimes \Cset((u))/
u\Cset[[u]],u\extd - u^{2}c)$ is a
quasi-isomorphism.

\item The chain map $\mathrm{Ch}\colon CC_\bullet(\A)
\otimes \Cset((u))\to \Omega^\bullet(X)\otimes \Cset((u))$
from the complex computing periodic cyclic homology
to the twisted de Rham complex $(\Omega^\bullet(X)
\otimes \Cset((u)),u\extd - u^{2}c)$ is a
quasi-isomorphism.
\end{enumerate}
\end{theorem}
\begin{pf*}{Proof.}
We shall prove 1 and 3 only, the proof of 2 is
identical and will be left to the reader.
Let us first choose a good open cover $\{U_i\}_{i\in I}$
of the compact manifold $X$.

To prove 1 above, we form the double complex
\begin{equation}
\label{eq:42}
\begin{matrix}
\hspace{-0.5 ex}\vdots & & \hspace{-0.5 ex}\vdots & & \hspace{-0.5 
ex}\vdots & & \\
\uparrow \scriptstyle{b} & &
\uparrow \scriptstyle{b} & &
\uparrow \scriptstyle{b} & &        \\
CC_2(\oplus \A_i) &
\stackrel{\d}{\rightarrow} &
CC_2(\oplus \A_{ij}) &
\stackrel{-\d}{\rightarrow} &
CC_2(\oplus \A_{ijk}) &
\stackrel{\d}{\rightarrow } & \cdots \\
\uparrow \scriptstyle{b} & &
\uparrow \scriptstyle{b} & &
\uparrow \scriptstyle{b} & &        \\
CC_1(\oplus \A_i) &
\stackrel{\d}{\rightarrow} &
CC_1(\oplus \A_{ij}) &
\stackrel{-\d}{\rightarrow} &
CC_1(\oplus \A_{ijk}) &
\stackrel{\d}{\rightarrow } & \cdots \\
\uparrow \scriptstyle{b} & &
\uparrow \scriptstyle{b} & &
\uparrow \scriptstyle{b} & &        \\
CC_0(\oplus \A_i) &
\stackrel{\d}{\rightarrow} &
CC_0(\oplus \A_{ij}) &
\stackrel{-\d}{\rightarrow} &
CC_0(\oplus \A_{ijk}) &
\stackrel{\d}{\rightarrow } & \cdots \\
\end{matrix}
\end{equation}
and notice that $\mathrm{Ch}$ defines a map
of double complexes between the double complex
\eqref{eq:42}
above and the following double complex.
\begin{equation}
\label{eq:43}
\begin{matrix}
\hspace{-0.5 ex}\vdots & & \hspace{-0.5 ex}\vdots & & \hspace{-0.5 
ex}\vdots & & \\
\uparrow \scriptstyle{0} & &
\uparrow \scriptstyle{0} & &
\uparrow \scriptstyle{0} & &        \\
\Omega^2(\coprod U_i) &
\stackrel{\d}{\rightarrow} &
\Omega^2(\coprod U_{ij}) &
\stackrel{-\d}{\rightarrow} &
\Omega^2(\coprod U_{ijk}) &
\stackrel{\d}{\rightarrow } & \cdots \\
\uparrow \scriptstyle{0} & &
\uparrow \scriptstyle{0} & &
\uparrow \scriptstyle{0} & &        \\
\Omega^1(\coprod U_i) &
\stackrel{\d}{\rightarrow} &
\Omega^1(\coprod U_{ij}) &
\stackrel{-\d}{\rightarrow} &
\Omega^1(\coprod U_{ijk}) &
\stackrel{\d}{\rightarrow } & \cdots \\
\uparrow \scriptstyle{0} & &
\uparrow \scriptstyle{0} & &
\uparrow \scriptstyle{0} & &        \\
\Omega^0(\coprod U_i) &
\stackrel{\d}{\rightarrow} &
\Omega^0(\coprod U_{ij}) &
\stackrel{-\d}{\rightarrow} &
\Omega^0(\coprod U_{ijk}) &
\stackrel{\d}{\rightarrow } & \cdots \\
\end{matrix}
\end{equation}
We now consider the two spectral sequences
associated to the double complex \eqref{eq:42}.
If we filter by rows
then it follows from Proposition \ref{prop: exact seq of complexes}
that the
$E_2$ term of the resulting spectral sequence is
just the Hochschild homology $HH_\bullet(\A)$ of
$\A$.  Now consider the spectral sequence associated
to \eqref{eq:42} obtained
by filtering by columns.  The map $\mathrm{Ch}$ of
double complexes described above yields a morphism
of spectral sequences.  Since the cover $\{U_i\}_{i\in I}$ is
good, each algebra $\oplus \A_{i_1\cdots i_p}$ is
isomorphic to the algebra $\oplus C^\infty(U_{i_1\cdots i_p})
\otimes \cL^1$.  In particular, in defining the map
$\mathrm{Ch}$ for each algebra $\oplus \A_{i_1\cdots i_p}$
we may take $\nabla$ to be the trivial connection $\extd$,
in which case the induced map $\mathrm{Ch}(\extd)\colon
HH_\bullet(\oplus \A_{i_1\cdots i_p})\to \Omega^\bullet(\coprod
U_{i_1\cdots i_p})$ coincides with
Connes' map~\eqref{eq: CHKR map},
which we know to be an isomorphism (Section \ref{sec:one b}).
By the homotopy invariance of the map $\mathrm{Ch}$
(Remark~\ref{rem: homotopy invariance of Hoch gen CHKR map})
the induced maps $\mathrm{Ch}(\nabla),
\mathrm{Ch}(\extd) \colon HH_\bullet(\oplus \A_{i_1\cdots i_p})
\to \Omega^\bullet(\coprod U_{i_1\cdots i_p})$ also coincide.
It follows then that the induced morphism of spectral
sequences is an isomorphism on the $E_1$ terms and hence
induces an isomorphism on the $E_\infty$ terms.  This
completes the proof of 1.

To prove 3, we first form the following double complex.
\begin{equation}
\label{eq:44}
\begin{array}{ccccccc}
\hspace{-3ex}\vdots & & \hspace{-4ex}\vdots & & \hspace{-4ex}\vdots & & 
\\
\hspace{1ex}\uparrow \scriptstyle{b+uB} & &
\uparrow \scriptstyle{b+uB} &  & \uparrow \scriptstyle{b+uB} & & \\
CC_1(\displaystyle{\bigoplus_i} \A_i)((u))
& \stackrel{\d}{\to} &
CC_1(\displaystyle{\bigoplus_{i<j}} \A_{ij})((u))
& \stackrel{-\d}{\to} &
CC_1(\displaystyle{\bigoplus_{i<j<k}} \A_{ijk})((u)) &
\stackrel{\d}{\to}  \cdots \\
\hspace{1ex}\uparrow \scriptstyle{b+uB} & &
\uparrow \scriptstyle{b+uB} & & \uparrow \scriptstyle{b+uB} & &  \\
CC_0(\displaystyle{\bigoplus_i} \A_i)((u)) &
\stackrel{\d}{\to} &
CC_0(\displaystyle{\bigoplus_{i<j}} \A_{ij})((u)) &
\stackrel{-\d}{\to} &
CC_0(\displaystyle{\bigoplus_{i<j<k}} \A_{ijk})((u)) &
\stackrel{\d}{\to}  \cdots  \\
\hspace{1ex}\uparrow \scriptstyle{b+uB} & &
\uparrow \scriptstyle{b+uB} & & \uparrow \scriptstyle{b+uB} & &  \\
CC_{-1}(\displaystyle{\bigoplus_i} \A_i)((u)) &
\stackrel{\d}{\to} &
CC_{-1}(\displaystyle{\bigoplus_{i<j}}\A_{ij})((u)) &
\stackrel{-\d}{\to} &
CC_{-1}(\displaystyle{\bigoplus_{i<j<k}}\A_{ijk})((u))
& \stackrel{\d}{\to}  \cdots \\
\hspace{1ex}\uparrow \scriptstyle{b+uB} & &
\uparrow \scriptstyle{b+uB} & & \uparrow \scriptstyle{b+uB} & &  \\
\hspace{-3ex}\vdots & & \hspace{-4ex}\vdots & & \hspace{-4ex}\vdots & &
\end{array}
\end{equation}
where $CC_k(\A)((u))$ denotes the elements of total
degree $k$ in $CC_k(\A)\otimes \Cset((u))$ as in \S~4.3.
Associated to this double complex are two spectral sequences
obtained by filtering by rows or by columns.  If we filter by
rows then it follows from
Proposition~\ref{prop: exact seq of complexes} above that the
$E_2$
term of the resulting spectral sequence is the periodic
cyclic homology of the algebra $\A$.  On the other hand, the map
$\mathrm{Ch}$ induces a map of bicomplexes from the
bicomplex~(\ref{eq:44})
to the bicomplex~(\ref{eq:20})
and hence a morphism of the spectral sequences
associated to these bicomplexes by filtering with respect to columns.
As in the proof of 1, it follows from
Proposition~\ref{prop: gen CHKR map indep of choice of connection}  and
Remark~\ref{rem: Hoch hom of C(X,L^1) = forms}
that this morphism of spectral
sequences is an isomorphism on the $E_1$ terms.  This proves the
3.
\end{pf*}

Recall that through the adjoint action the projective unitary group
$PU$ acts on the Schatten ideals $\cL^p$, $p\geq 1$.
Therefore for any $p\geq 1$ we can form the associated bundle
$\cL^p(P)$ and the Fr\'{e}chet algebra of smooth sections
$C^\infty(\cL^p(P))$.
The $K$-theory of the algebra $C^\infty(\cL^p(P))$ is the
twisted $K$-theory $K_i(C(\K(P)))$.  We can, with very little
extra effort, compute the periodic cyclic homology of the
algebra $C^\infty(\cL^p(P))$ using Theorem~\ref{thm: main result}.
The inclusion $\cL^1\hookrightarrow \cL^p$ induces an inclusion
$\cL^1(P) \hookrightarrow \cL^p(P)$.  Associated to the two algebras
$C^\infty(\cL^1(P))$ and $C^\infty(\cL^p(P))$ and
the good cover $\{U_i\}_{i\in I}$ of $X$ we have double complexes
of the form \eqref{eq:44}.
The inclusion $\cL^1(P)\hookrightarrow \cL^p(P)$ induces a morphism
of these double complexes.
Recall that the inclusion
$\cL^1\hookrightarrow \cL^p$ induces an isomorphism on periodic
cyclic homology.  More generally we have (\cite{Cu1} Corollary 17.2)
if $U$ is an open subset of $X$ then
the inclusion $C^\infty(U)\otimes \cL^1 \hookrightarrow
C^\infty(U)\otimes \cL^p$ induces an isomorphism on periodic
cyclic homology.
Consider the double complexes~(\ref{eq:44})
associated
to the algebras $C^\infty(\cL^1(P))$ and $C^\infty(\cL^p(P))$.  If we
filter
the respective total complexes by columns then the morphism
$C^\infty(\cL^1(P))\hookrightarrow C^\infty(\cL^p(P))$ induces an
isomorphism on the $E_1$ terms of the resulting spectral sequences.
It follows then that $HP_i(C^\infty(\cL^1(P)))$ and
$HP_i(C^\infty(\cL^p(P)))$
are isomorphic.
\begin{corollary}
For $p\geq 1$ the inclusion
$C^\infty(\cL^1(P))\hookrightarrow C^\infty(\cL^p(P))$ induces
an isomorphism on periodic cyclic homology.
\end{corollary}

%-----------------------------------------------------------------------

\section{The Twisted Chern Character}

In this section we relate the Connes-Chern character
$\mathrm{ch}\colon K_\bullet(\A)\to HP_\bullet(\A)$
to the twisted Chern character developed in
earlier work \cite{BCMMS}.  Notice that we can view
$\mathrm{ch}$ as a homomorphism $\mathrm{ch}_P\colon
K^\bullet(X,P)\to H^\bullet(X,c(P))$ once we notice
the identifications $K_\bullet(\A)\cong K^\bullet(X,P)$
and $HP_{\bullet}(\A)\cong H^\bullet(X,c(P))$
where the first isomorphism is established in
\cite{Ros} and the second isomorphism
is Theorem~\ref{thm: main result}.

We first observe
that the Connes-Chern character $\mathrm{ch}$, or
alternatively the homomorphism $\mathrm{ch}_P \colon
K^\bullet(X,P)\to H^\bullet(X,c(P))$, becomes an
isomorphism after tensoring with the complex numbers.
Recall that in general one has no assurances that the Connes-Chern map
$\text{ch}\colon K_i(\A) \to HP_i(\A)$ will be an
isomorphism after tensoring with $\Cset$
for arbitrary algebras $\A$ --- for example
$K_i(C(X)) = K^i(X)$ but $HP_i(C(X)) = 0$ where
$C(X)$ is the algebra of continuous complex
valued functions on a compact space $X$.  However, in our situation,
we can prove the following.

\begin{proposition}\label{connes-chern iso}
For the algebra $\A = C^\infty(\cL^1(P))$,
the Connes-Chern character $\;\mathrm{ch}\colon K_i(\A)
\to HP_i(\A)$ becomes an isomorphism after tensoring
with $\Cset$, for $i=0,1$.
\end{proposition}

\begin{pf*}{Proof.}
We prove this Proposition by making repeated use of the
Mayer-Vietoris exact sequence (note first that the construction
of the Mayer-Vietoris sequence given in \cite{Black}, pages 219--220,
goes across to the periodic cyclic theory $HP$).
First of all choose a cover $\{B_{\a}\}_{\a\in A}$
of $X$ by closed, contractible sets $B_{\a}$ (this
can be done using a minor modification of the usual
construction of a good open cover by geodesically convex
balls).  Consider the Mayer-Vietoris sequence
for the pair of closed sets $\{B_{\a},B_{\b}\}$.  Let
us write $\A_{B_{\a}}$ and $\A_{B_{\a\b}}$ for the algebras
$C^\infty(\cL^1(P)|_{B_{\a}})$ and
$C^\infty(\cL^1(P)|_{B_{\a\b}})$ respectively,
where $B_{\a\b} = B_{\a}\cap B_{\b}$.
Then the Mayer-Vietoris exact sequence for $\A|_{B_{\a}\cup B_{\b}}$ in
$K$-theory is
\begin{multline}
\label{eq:15}
\rightarrow K_i(\A|_{B_{\a}\cup B_{\b}}) \rightarrow
K_i(\A|_{B_{\a}})\oplus K_i(\A|_{B_{\b}}) \rightarrow
K_i(\A|_{B_{\a\b}}) \rightarrow K_{i+1}(\A|_{B_{\a}\cup B_{\b}})
\rightarrow
\end{multline}
and there is a similar exact sequence for the periodic cyclic
theory.  From results of Nistor \cite{Nis}, we
see that the Connes-Chern map $\text{ch}$
furnishes us with
a morphism from the exact sequence \eqref{eq:15}
to the corresponding sequence for $HP$.  We know that
$\text{ch}\colon K_i(\A|_{B_{\a}})\to HP_i(\A|_{B_{\a}})$
and $\text{ch}\colon K_i(\A|_{B_{\a\b}})\to
HP_i(\A|_{B_{\a\b}})$ become isomorphisms after tensoring with
$\Cset$ for all $\a$ and $\b$.  It follows therefore by the
$5$-lemma that $\text{ch}\colon K_i(\A|_{B_{\a}\cup B_{\b}})
\to HP_i(\A|_{B_{\a}\cup B_{\b}})$ is an isomorphism
after tensoring with $\Cset$.  By induction we see that
$\text{ch}\colon K_i(\A)\otimes \Cset \to HP_i(\A)\otimes \Cset$
is an isomorphism.
\end{pf*}

Upon identifying $K_i(\A)$ with $K^i(X,P)$ (see \cite{BCMMS})
and $HP_i(\A)$ with $H^i(X,c(P))$ (see Theorem \ref{main}), we denote 
the
Connes-Chern character $\;\mathrm{ch}\colon K_i(\A)
\to HP_i(\A)$ by $\widetilde{\mathrm{ch}}_P\colon
K^i(X,P)\to H^i(X,c(P))$.
The following Proposition is easy to establish.

\begin{proposition}
The Connes-Chern character $\widetilde{\mathrm{ch}}_P\colon
K^i(X,P)\to H^i(X,c(P))$ satisfies the
following properties:
\begin{enumerate}
\item If the bundle $P$ is trivial (say
$P$ is trivialised by a section $s$)
then $\widetilde{\mathrm{ch}}_P$ reduces to the
ordinary Chern character $\mathrm{ch}
\colon K^i(X)\to H^i(X)$ in the sense
that we have the commutative diagram
$$
\xymatrix{
K^i(X,P) \ar[r]^{\widetilde{\mathrm{ch}}_P}
\ar[d]_-{\cong} & H^i(X,c(P))
\ar[d]^-{e^{u s^*f}}       \\
K^i(X) \ar[r]_-{\mathrm{ch}}
& H^i(X)                          }
$$
\item $\widetilde{\mathrm{ch}}_P\colon K^i(X,P)
\to H^i(X,c(P))$ is natural with respect
to maps $h\colon Y\to X$ in the sense
that we have the following
commutative diagram
$$
\xymatrix{
K^i(X,P) \ar[r]^-{\widetilde{\mathrm{ch}}_P}
\ar[d]_-{h^*} & H^i(X,c(P))
\ar[d]^-{h^*}                   \\
K^i(Y,h^*P) \ar[r]^-{\widetilde{\mathrm{ch}}_{h^*P}}
& H^i(Y,h^*c(P)). }
$$
\end{enumerate}
\end{proposition}
In $(1)$ above, the de Rham cocycle $c(P)$ is trivialised
by $s$ since we have $c(P) = \extd s^*f$, where $f$ is
the $2$-form~\eqref{eq:16}.

In \cite{BCMMS}, a geometric,  twisted Chern character was defined,
\begin{equation}
\label{eq:104}
\mathrm{ch}_P : K^\bullet (X, P) \to H^\bullet (X, c(P)),
\end{equation}
using connections and curvature as well as a curving for the principal 
$PU$ bundle $P$,
which has the same functorial properties as the Connes-Chern character,
namely that is natural with respect to maps and that it reduces to the
ordinary
Chern character when $P$ is the trivial $PU$ bundle.  The twisted Chern
character
also satisfies other functorial properties, such as being compatible
with the $K^0(X)$
module structure of $K^\bullet (X, P) $ and the $H^{even} (X)$ module
structure of
$H^\bullet (X, c(P))$, cf. \cite{MS}. The analogous properties for the
Connes-Chern character
can be deduced from the following basic result that identifies the two
homomorphisms.

\begin{proposition}
For the algebra $\A = C^\infty(\cL^1(P))$,
the Connes-Chern isomorphism $\mathrm{ch}\colon K_i(\A)\otimes \Cset
\to HP_i(\A) \otimes \Cset $  coincides with the twisted Chern character
ismorphism  $\mathrm{ch}_P \colon K^i(X, P)\otimes \Cset
\to H^i(X, c(P)) \otimes \Cset$ for $i=0,1$.
\end{proposition}

\begin{pf*}{Proof.}
We follow the strategy of Proposition~\ref{connes-chern iso}.
Let $\{B_{\a}\}_{\a\in A}$
be a (finite) cover of the compact manifold $X$ by
closed balls $B_{\a}$ such that $\cL^1(P)$ is
trivialised over $B_{\a}$.  Consider the
Mayer-Vietoris exact sequences in $K$-theory and periodic cyclic
homology for the algebra $\A|_{B_{\a}}$. A choice of local section
$s_\alpha : B_\alpha \to P\big|_{B_\alpha}=P_\alpha$, trivializes
$\A|_{B_{\a}}
\cong
C^\infty (B_\alpha) \otimes \cL^1$.  Using this local section, we will
have
the commutative
diagram

\begin{equation}\label{cc-iso}
\begin{CD}
K_\bullet( C^\infty (B_\alpha) \otimes \cL^1)  @>\cong>> K_\bullet(
C^\infty (B_\alpha))@>\cong>>
K^\bullet(B_\alpha) @>\cong>>  K^\bullet(B_\alpha, P_\alpha) \\
@V{\rm ch} VV @V{\rm ch}VV @V{\rm ch}VV @V{{\rm ch}_{P_\alpha}}VV
                        \\
HP_\bullet( C^\infty (B_\alpha) \otimes \cL^1)  @>\cong>> HP_\bullet(
C^\infty (B_\alpha)) @>\cong>> H^\bullet(B_\alpha) @>\cong>>
H^\bullet(B_\alpha, c(P_\alpha))
\end{CD}
\end{equation}

where the horizontal arrows are all isomorphisms, and the vertical
arrows
are isomorphisms after tensoring with $\Cset$. The first and second
vertical arrows
starting from the left are the Connes-Chern isomorphism. The third
vertical arrow is
the usual Chern character, which is well known can be identified with
the Connes-Chern
character on smooth functions on a manifold. The last isomorphism uses
the local section
again, where $s_\alpha^*f = f_\alpha$ and $s_\alpha^*c =  c_\alpha =
df_\alpha$,
and the last vertical arrow is the twisted Chern character.
The first upper horizontal isomorphism is essentially Morita 
equivalence.
More precisely,
one uses fact that $C^\infty (B_\alpha) \otimes \cL^1$ is a smooth
subalgebra of $C(B_\alpha) \otimes \K$, so that the inclusion map
induces an isomorphism in $K$-theory. Then one uses the invariance of
$K$-theory upon tensoring with the compact operators $\K$.
The second upper horizontal isomorphism is the Serre-Swan
isomorphism, the third is the fact the the local section trivializes
the $PU$ bundle over $B_\alpha$
and we use Morita invariance again.
The first lower horizontal isomorphism is Morita equivalence, the
second is the Connes-Hochschild-Kostant-Rosenberg isomorphism, the
third isomorphism is the fact that the twisted deRham differential is
obtained by
conjugating the standard one by $e^{uf_\alpha}$. In summary, a local
section
$s_\alpha : B_\alpha \to P\big|_{B_\alpha}=P_\alpha$ gives rise to the
commutative diagram,

$$
\begin{CD}
K_\bullet( \A|_{B_{\a}})  @>\cong>>   K^\bullet(B_\alpha, P_\alpha) \\
@V{\rm ch} VV @V{{\rm ch}_{P_\alpha}}VV                        \\
HP_\bullet( \A|_{B_{\a}})  @>\cong>> H^\bullet(B_\alpha, c(P_\alpha))
\end{CD}
$$
where the upper horizontal arrow is the isomorphism in \cite{BCMMS}
and is the result of the composition of the upper horizontal
isomorphisms in
equation \eqref{cc-iso}:
the lower horizontal arrow is the generalized
Connes-Hochschild-Kostant-Rosenberg isomorphism that is established in
this paper,
and is the result of the composition of the lower horizontal
isomorphisms in
equation \eqref{cc-iso},
We next use the Mayer-Vietoris sequence, the 5-lemma  and induction as
in Proposition \ref{connes-chern iso} to establish the proposition.
\end{pf*}
%-----------------------------------------------------------------------

%-----------------------------------------------------------------------

\end{document}